\documentclass[12pt]{article}         %[dds-c.tex 29.10.09]
\usepackage{amsfonts,amssymb} %mlb,
\makeatletter
\def\bdraft{\pagestyle{myheadings} %[draft + upper page numbers]
           \textheight=10.5truein \textwidth=7.5truein \parindent=8pt
           \voffset=-1truein \topmargin=30pt \headheight=10pt \headsep=3pt
           \ifcase \@ptsize \hoffset=-1.5truein \or \hoffset=-1.35truein
                        \or \hoffset=-1.15truein \fi}
\def\quality{\textheight=240mm \textwidth=160mm \topmargin=0Truein
             \ifcase \@ptsize \hoffset=-23mm
                     \or \hoffset=-20mm \or \hoffset=-15mm \fi}
\def\USquality{\textheight=250mm \textwidth=170mm \topmargin=0Truein
               \voffset=-1truein
               \ifcase \@ptsize \hoffset=-23mm
                       \or \hoffset=-20mm \or \hoffset=-15mm \fi}
\makeatother \USquality %\bdraft
\def\beq#1#2{\begin{equation} \label{#1} #2 \end{equation}}
\def\bea#1{\begin{eqnarray*} #1 \end{eqnarray*}} \def\a{\!\!\!&\!\!\!\!&}
\def\n{\noindent}   \def\map{T}   
    
\def\IR{{\mathbb{R}}}  \def\IZ{{\mathbb{Z}}}  
\def\toas#1{\stackrel{#1}{\longrightarrow}}
\def\ep{\varepsilon}  \def\phi{\varphi}   
  \def\mod1{\,({\rm mod\ } 1)\,}
\def\t{\tilde} \def\den{\rho}   

\def\function#1{\left\{\!\!\!\begin{array}{ll} #1 \end{array} \right.}
\def\proof{\smallskip \noindent {\bf Proof. \ }}       %start of proof
\def\blanksquare{\,\,\,$\sqcup\!\!\!\!\sqcap$}         %blank  square
\def\qed{\hfill\blanksquare\linebreak\smallskip\par}   %end of proof

\def\thname{Theorem}     \def\lmname{Lemma}      \def\prname{Proposition}
\def\dfname{Definition}  \def\crname{Corollary}  \def\rmname{Remark}

\newtheorem{theorem}{\thname}[section]   %Numbering: Theorem--Other section
\newtheorem{lemma}{\lmname}[section]     %{lemma}[theorem]{Lemma}   section
\newtheorem{corollary}[lemma]{\crname}   %lemma

\newtheorem{dftn}{\dfname}[section]
\newtheorem{rmrk}[lemma]{\rmname}
\newenvironment{remark}{\begin{rmrk}\rm}{\end{rmrk}}     %lemma

             \catcode`@=11 \@addtoreset{equation}{section} \catcode`@=12

\newcommand\mlbscale{1pt} %to change: \renewcommand\mlbscale{1.3pt}
\newif\iffigs\figstrue %\newif\iffigs\figsfalse -- to fake figures
\def\bline(#1,#2)(#3,#4)(#5){\put(#1,#2){\line(#3,#4){#5}}}  %straight line

\def\bfig(#1,#2)#3#4{\begin{figure} \begin{center}
    \framebox{\setlength{\unitlength}{\mlbscale}
       \iffigs \begin{picture}(#1,#2) #3 \end{picture}
       \else \begin{picture}(60,10)(0,0)
                   \put(0,0){\framebox(60,10){Figure}} \end{picture} \fi}
    \end{center} \caption{#4} \end{figure}}
\def\Bfig(#1,#2)#3#4{\begin{figure} \begin{center}
    \setlength{\unitlength}{\mlbscale}
       \iffigs \begin{picture}(#1,#2) #3 \end{picture}
       \else \begin{picture}(60,10)(0,0)
                   \put(0,0){\framebox(60,10){Figure}} \end{picture} \fi
    \end{center} \caption{#4} \end{figure}}
\def\bpic(#1,#2)#3{\setlength{\unitlength}{\mlbscale}
    \begin{picture}(#1,#2) #3 \end{picture}}
   \def\cV{{\cal V}}
\def\cN{{\cal N}} \def\cNs{\cN_s}  \def\cNw{\cN_w}  
\def\cB{{\cal B}} \def\cC{{\cal C}}
\def\*#1{#1^*}    \def\0#1{\breve#1}  \def\2#1{\acute#1}
\def\G{\Delta}  % {G}

\def\?#1{} % comments

\begin{document}%%%-------------------------------------------------%%%%
\title{Metric properties of discrete time \\
       exclusion type processes in continuum}
\author{Michael Blank\thanks{
        Russian Academy of Sci., Inst. for
        Information Transm. Problems,
        and Laboratoire Cassiopee UMR6202, CNRS, France; ~
        e-mail: blank@iitp.ru}
        %\thanks{This research has been partially supported
        %        by Russian Foundation for Fundamental Research,
        %        Program ONIT, and French Ministry of Education grants.}
       }
\date{October 29, 2009} %\today}%{April 26, 2009} %
\maketitle

\begin{abstract}%
A new class of exclusion type processes acting in continuum with
synchronous updating is introduced and studied. Ergodic averages
of particle velocities are obtained and their connections to other
statistical quantities, in particular to the particle density (the
so called Fundamental Diagram) is analyzed rigorously. The main
technical tool is a ``dynamical'' coupling applied in a
nonstandard fashion: we do not prove the existence of the
successful coupling (which even might not hold) but instead use
its presence/absence as an important diagnostic tool. Despite that
this approach cannot be applied to lattice systems directly, it
allows to obtain new results for the lattice systems embedding
them to the systems in continuum. Applications to the traffic
flows modelling are discussed as well.
\end{abstract}%

%%%%%%%%%%%%%%%%%%%%%%%%%%%%%%%%
\section{Introduction}\label{s:intro}%{Introduction}

The classical simple exclusion process is a Markov process that
describes nearest-neighbor random walks of a collection of
particles on the one-dimensional infinite\footnote{or finite with
   periodic boundary conditions} %
integer lattice. Particles interact through the hard core
exclusion rule, which means that at most one particle is allowed
at each site. This seemingly very particular process introduced
first in 1970 by Frank Spitzer \cite{Sp} appears naturally in a
very broad list of scientific fields starting from various models
of traffic flows \cite{NS,GG,ERS,Bl-erg,Bl-hys}, molecular motors
and protein synthesis in biology(see e.g. \cite{SZL}), surface
growth or percolation processes in physics (see \cite{Pe,BFS} for
a review), and up to the analysis of Young diagrams in
Representation Theory \cite{CMOS}.

Qualitatively from the point of view of the order of particle
interactions there are two principally different types of exclusion
processes: with synchronous and asynchronous updating rules. In the
latter case at each moment of time a.s. at most one particle may
move and hence only a single interaction may take place. This is the
main model considered in the mathematical literature (see e.g.
\cite{Lig,Thor,Sp,Num} for a general account and \cite{An,EFM,Ros}
for recent results), and indeed, the assumption about the
asynchronous updating is quite natural in the continuous time
setting. The synchronous updating means that {\em all} particles are
trying to move simultaneously and hence an arbitrary large (and even
infinite) number of interactions may occur at the same time. This
makes the analysis of the synchronous updating case much more
difficult, but this is what happens in the discrete time
case.\footnote{if one do not consider some
   ``artificial'' updating rules like a sequential or random updating.} %
This case is much less studied, but still there are a few results
describing ergodic properties of such processes
\cite{BF,Bl-erg,Bl-hys,BP,ERS,GG,NS}.

Our aim is to introduce and study the synchronous updating version
of the exclusion process in continuum. Note that recently some
other interacting particle processes were generalized from
lattice to continuum case (see e.g. \cite{Pe,FKO}).

A configuration $x:=\{x_i\}_{i\in\IZ}$ is a bi-infinite sequence
of real numbers $x_i\in\IR$ interpreted as centers of particles
represented by balls of radius $r\ge0$ (see Fig.~\ref{f:tasep-c})
and ordered with respect to their positions (i.e. $\dots\le x_{-1}
\le x_0 \le x_1 \le \dots$). To emphasize the dependence on the
radius $r\ge0$ we shall use the notation $x(r)$. We say that a
configuration $x(r)$ is {\em admissible} if
$$x_i(r)+r\le x_{i+1}(r)-r~~\forall i\in\IZ$$
(the corresponding balls may only touch each other) and
denote by $X(r)$ the space of admissible configurations.

%%%%%%%%%%%%%%%%%%%%%%%%%%%%%%%%%%%%%%%%%%%%%%%%%%%%%%%%%%%%%%%%%%
%% TASEP in continuum
\def\particle#1{
     \put(0,0){\circle{30}} \bline(0,0)(0,-1)(25) \put(-3,-32){$x_{#1}$}
     \bline(0,20)(1,0)(15) \put(5,22){$r$}
     \put(0,5){\vector(1,0){40}} \put(20,9){$v_{#1}$}}
\Bfig(150,30)
      {\footnotesize{
       \thicklines
       \bline(0,10)(1,0)(150)
       \put(30,10){\particle{i}} \put(110,10){\particle{i+1}}
       \thinlines \bline(46,-18)(1,0)(48) \put(65,-13){$\Delta_i$}
       %\bline(94,-10)(0,1)(40)
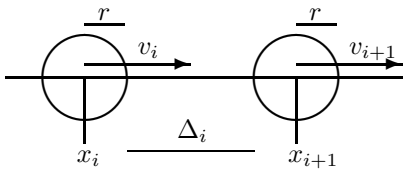
      }      } {TASEP in continuum. \label{f:tasep-c}}
%%%%%%%%%%%%%%%%%%%%%%%%%%%%%%%%%%%%%%%%%%%%%%%%%%%%%%%%%%%%%%%%%

The dynamics will be defined as follows. We assume given a
collection of (possibly random) values $\{v_i^t\}_{i,t}$, where
$i,t\in\IZ$ and $t\ge0$; conditions on this collection will be
given shortly. For a trivial configuration consisting of a single
particle located at time $t\ge0$ at $x_0^t\in\IR$ (i.e.
$x^t\equiv\{x_0^t\}$) the dynamics is defined as %
$$x_0^{t+1}:=x_0^t+v_0^t,$$
and thus $v_0^t$ may be interpreted as a local velocity at time
$t$, i.e. this is simply a random walk on $\IR$. To generalize
this trivial setting for an infinite configuration $x(r)\in X(r)$
we again interpret a (be-infinite on $i\in\IZ$) sequence
$\{v_i^t\}_{i,t}$ as {\em local velocities} for particles in
$x^t(r)$ performing random walks conditioned to the order
preservation and the hard core exclusion rule.

To simplify presentation we restrict ourselves here to the case of
nonnegative local velocities postponing the discussion of the
general case when the local velocities take both positive and
negative signs to Section~\ref{s:vel-2signs}. The point is that
the formulations in the latter case are becoming much more
involved, but the results and arguments work with only very slight
changes.

Since only nonnegative local velocities are considered the hard
core exclusion rule means that the admissibility condition breaks
down for the $i$-th particle at time $t\in\IZ_+$ if and only if
the inequality %
$$x^t_i(r)+v_i^t + r \le x^t_{i+1}(r) - r$$ %
does not hold. If this happens we say that there is a {\em
conflict} between the particles $i$ and $i+1$, and to resolve it
one applies a {\em normalizing} construction
$$v_i^t\to\cN(v_i^t,x^t(r)).$$ After the normalization the positions
of particles are calculated according to the rule
$$x_i^{t+1}(r):=x_i^t(r)+\cN(v_i^t,x^t(r))~~\forall i.$$

In what follows we always assume\footnote{This formulation allows to
   consider velocities of both signs which we shall do in
   Section~\ref{s:vel-2signs} and simply means that the normalized
   velocity has the same direction as the original one and cannot
   exceed it on modulus.} %
that $\forall i,t~~ \cN(v_i^t,x^t(r))\in[0,v_i^t]$ (to simplify notation by
   the segment $[a,b]$ we mean $[\min(a,b),\max(a,b)]$) %
and consider only {\em nonanticipating} normalizations\footnote{In
   Section~\ref{s:gen} we shall show that the violation of this
   condition makes the system to be not well posed.} %
satisfying the condition that in the case of the conflict of the
$i$-th particle with the $j$-th one\footnote{For nonnegative
   velocities $j\equiv i+1$, but in general $j\in\{i-1,i+1\}$.} %
at time $t$ the position of the $i$-th particle at the next moment
of time $x_i^{t+1}(r)\in[x_i^t(r),x_j^t(r)]$.

The normalization may be done in a number of ways and we restrict
ourselves to two extreme constructions. The first of them we call
{\em strong normalization} (notation $\cNs(\cdot,\cdot)$) and
according to the name we reject (nullify) the velocity leading to
the conflict. The second construction we call {\em weak
normalization} (notation $\cNw(\cdot,\cdot)$) and in this case we
modify the conflicting velocity to allow the particle to move as
far as possible. In terms of {\em gaps}\/ %
$$\G_i(x^t(r))\equiv\G_i^t:=x_{i+1}^t(r)-x_{i}^t(r)-2r$$
between particles in the configuration $x^t$ the normalization
procedures are written as follows: %
$$ \cNs(v_i^t,x^t(r)):=\function{v_i^t &\mbox{if } v_i^t\le \G_i^t  \\
                         0     &\mbox{otherwise }}, \qquad
   \cNw(v_i^t,x^t(r)):=\function{v_i^t &\mbox{if } v_i^t\le \G_i^t  \\
                         \G_i^t &\mbox{otherwise }} .$$
Fig.~\ref{f:normalization} demonstrates possible positions of
particles at two consecutive moments of time $t$ and $t+1$ for the
cases of weak (a-c) and strong (a'-c') normalizations. Despite
appearances these two normalization procedures lead to a very
different limit behavior of the corresponding particle systems.
The simplest example (existing even in the continuous time case)
is the situation when $v_i^t\equiv{v}~\forall i,t$ and the gaps
between particles in $x$ are smaller than $v$. Then under the
strong normalization no motion is allowed, while the weak
normalization leads to the well defined motion -- the exchange of
gaps between particles.
Other normalization procedures together with more general
assumptions about the dynamics will be discussed in
Section~\ref{s:gen}.

%%%%%%%%%%%%%%%%%%%%%%%%%%%%%%%%%%%%%%%%%%%%%%%%%%%%%%%%%%%%%%%%%%
%% Normalization
\Bfig(290,120)
      {\footnotesize{
       %%%%% Weak %%%%%
       %% picture (a)
       \put(0,100){\bpic(100,30){\put(-20,5){(a)}
          \bline(0,0)(1,0)(100) \bline(0,10)(1,0)(100)
          \put(10,10){\circle{5}} \put(70,10){\circle{5}}
          \put(10,15){\vector(1,0){40}} \put(70,15){\vector(1,0){20}}
          \put(50,0){\circle{5}} \put(90,0){\circle{5}}
          \put(105,7){$t$} \put(105,-3){$t+1$}
          \put(10,18){$i$}  \put(60,18){$i+1$}
          \put(47,-12){$i$} \put(80,-12){$i+1$}
          }}
       %% picture (b)
       \put(0,50){\bpic(100,30){\put(-20,5){(b)}
          \bline(0,0)(1,0)(100) \bline(0,10)(1,0)(100)
          \put(10,10){\circle{5}} \put(70,10){\circle{5}}
          \put(10,15){\vector(1,0){70}} \put(70,13){\vector(1,0){20}}
          \put(65,0){\circle{5}} \put(90,0){\circle{5}}
          \put(105,7){$t$} \put(105,-3){$t+1$}
          \put(10,18){$i$}  \put(60,18){$i+1$}
          \put(63,-12){$i$} \put(80,-12){$i+1$}
          }}
       %% picture (c)
       \put(0,0){\bpic(100,30){\put(-20,5){(c)}
          \bline(0,0)(1,0)(100) \bline(0,10)(1,0)(100)
          \put(10,10){\circle{5}} \put(70,10){\circle{5}}
          \put(10,15){\vector(1,0){50}} \put(70,13){\vector(-1,0){35}}
          \put(45,0){\circle{5}} \put(50,0){\circle{5}}
          \put(105,7){$t$} \put(105,-3){$t+1$}
          \put(10,18){$i$}  \put(60,18){$i+1$}
          %\put(47,-12){$i$} \put(80,-12){$i+1$}
          }}
       %%%%% Strong %%%%%
       %% picture (a')
       \put(170,100){\bpic(100,30){\put(-20,5){(a')}
          \bline(0,0)(1,0)(100) \bline(0,10)(1,0)(100)
          \put(10,10){\circle{5}} \put(70,10){\circle{5}}
          \put(10,15){\vector(1,0){40}} \put(70,15){\vector(1,0){20}}
          \put(50,0){\circle{5}} \put(90,0){\circle{5}}
          \put(105,7){$t$} \put(105,-3){$t+1$}
          \put(10,18){$i$}  \put(60,18){$i+1$}
          \put(47,-12){$i$} \put(80,-12){$i+1$}
          }}
       %% picture (b')
       \put(170,50){\bpic(100,30){\put(-20,5){(b')}
          \bline(0,0)(1,0)(100) \bline(0,10)(1,0)(100)
          \put(10,10){\circle{5}} \put(70,10){\circle{5}}
          \put(10,15){\vector(1,0){70}} \put(70,13){\vector(1,0){20}}
          \put(10,0){\circle{5}} \put(90,0){\circle{5}}
          \put(105,7){$t$} \put(105,-3){$t+1$}
          \put(10,18){$i$}  \put(60,18){$i+1$}
          \put(10,-12){$i$} \put(80,-12){$i+1$}
          }}
       %% picture (c')
       \put(170,0){\bpic(100,30){\put(-20,5){(c')}
          \bline(0,0)(1,0)(100) \bline(0,10)(1,0)(100)
          \put(10,10){\circle{5}} \put(70,10){\circle{5}}
          \put(10,15){\vector(1,0){50}} \put(70,13){\vector(-1,0){35}}
          \put(10,0){\circle{5}} \put(70,0){\circle{5}}
          \put(105,7){$t$} \put(105,-3){$t+1$}
          \put(10,18){$i$}  \put(60,18){$i+1$}
          \put(10,-12){$i$} \put(60,-12){$i+1$}
          }}
      }}
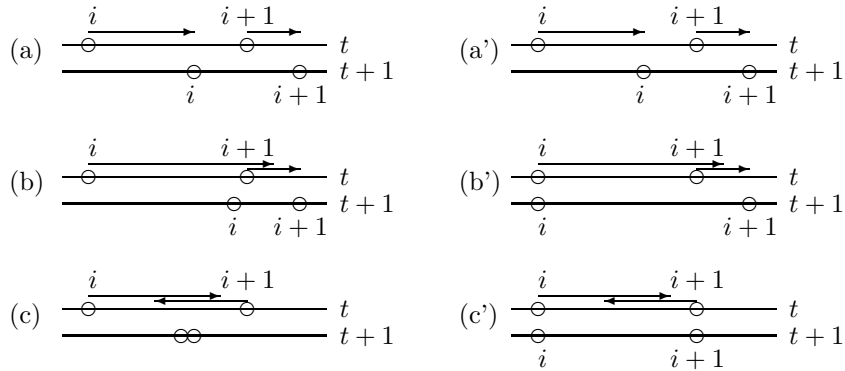
{Positions of particles at time $t,t+1$ in cases of weak (a-c) and
strong (a'-c') normalizations. Local particle velocities are shown
by vectors. The cases (c,c') correspond to negative velocities and
will be discussed in Section~\ref{s:vel-2signs}. \label{f:normalization}}
%%%%%%%%%%%%%%%%%%%%%%%%%%%%%%%%%%%%%%%%%%%%%%%%%%%%%%%%%%%%%%%%%

Observe that any two particle configurations $x(r),~\2x(\2r)$
having the same sequence of gaps $\G:=\{\G_i\}$ may be transformed
to each other by a one-to-one map %
\beq{e:r->r'}{\2x_i(\2r)=\phi(x_i(r)):=x_i(r)-2i(r-\2r)~~\forall i\in\IZ.} %
Since the normalization procedures that we consider depend only on
the gaps between particles it is enough to study the case $r=0$. On
the other hand, if $r=1/2, ~x_i^0(r)\in\IZ~\forall i\in\IZ$ and
$v_i^t\in\IZ~\forall i\in\IZ, t\ge0$ then $x_i^t(r)\in\IZ~\forall
i\in\IZ, t\ge0$ which means that we get a lattice particle system.
Thus our results lead to a completely new approach to the analysis of
lattice systems as well. Note however that in the case $r=0$ an arbitrary
number of particles may share the same spatial position which is
prohibited in the lattice case.

Due to the observation above we shall study in detail only the case $r=0$
since the corresponding results for any $r>0$ are readily available
through the transformation (\ref{e:r->r'}), see e.g. specific
calculations for densities and velocities in
Lemmas \ref{l:den-r}, \ref{l:V-r} and Corollaries~\ref{c:V-r-w}, \ref{c:V-r-s}.

To simplify notation we shall use the convention
$x(r)\equiv x^0(r), ~x\equiv x^0(0)$ and similarly $X\equiv X(0)$.

Of course, without some specific assumptions on the structure of
local velocities $\{v_i^t\}_{i,t}$ no interesting results are
possible. We assume that %
$v_i^t\in[0,v]~~\forall i\in\IZ,~t\in\IZ_0:=\IZ_+\cup\{0\}$ and
one of the following seemingly opposite assumptions holds: %
\begin{itemize}
\item[(a)] $v_i^t\equiv v_0^t~~\forall i\in\IZ,~t\in\IZ_0$ and
    $\exists~\bar{v}(\gamma):=\lim\limits_{t\to\infty}\frac1t
                      \sum\limits_{s=0}^{t-1}\min(v_0^s,\gamma)
     ~~\forall\gamma>0$ ~~(a.s.); %
\item[(b)] $\{v_i^t\}$ are i.i.d. (both in $i$ and $t$) random
    variables.
\end{itemize}

Note that the intersection between the sets of local velocities
satisfying the assumptions (a) and (b) contains an important case
of pure deterministic velocities: $v_i^t\equiv v~~\forall
i\in\IZ,~t\in\IZ_0$. As we shall show properties of systems with
local velocities satisfying to the assumption (a) are close to the
pure deterministic setting. Therefore we refer to the setting (a)
as {\em deterministic}\footnote{In this case
   $v_0^t$ might be a trajectory of a deterministic chaotic map
   $f:[0,1]\to[0,1]$, e.g. $v_0^{t+1}:=vf^t(v_0^t/v)$, as well as
   a realization of a true random Markov chain).} %
and to the setting (b) as {\em random}.

It is of interest that in the seemingly simplest purely
deterministic setting $v_i^t\equiv v~~\forall i\in\IZ,~t\in\IZ_0$
the behavior of the corresponding deterministic dynamical system
describing the dynamics of particle configurations is far from
being trivial. In Section~\ref{s:entropy} we prove that this
system is chaotic in the sense that its topological entropy is
positive (and even infinite).

To emphasize that under dynamics no creation or annihilation of
particles may take place this sort of systems is called {\em
diffusive driven systems} (DDS) instead of a more general object
-- {\em interacting particle systems} (IPS).

The main technical tool in our analysis is a (somewhat unusual)
``dynamical'' coupling construction. Despite that various
couplings are widely used in the analysis of IPS, applications of
our approach is very different from conventional. In particular,
we do not prove the existence of the so called successful coupling
(which even might not hold) but instead use its presence/absence
as an important diagnostic tool. Remark also that typically one
uses the coupling argument to prove the uniqueness of the
invariant measure and to derive later other results from this
fact. In our case there might be a very large number of ergodic
invariant measures or no invariant measures at all (recall the
trivial example of a single particle performing a skewed random
walk). The latter example indicates that there is another
important statistical quantity -- average particles velocity that
can be computed at least in this case. (See e.g. \cite{BBM} for a
discussion of the average velocity in the context of Queueing
Networks.) The dynamical coupling will be used directly to find
connections between the average particle velocities and other
statistical features of the systems under consideration, in
particular with the corresponding particle densities.

It is worth note that all approaches used to study lattice
versions of DDS are heavily based on the combinatorial structure
of particle configurations. This structure has no counterparts in
the continuum setting under consideration. In particular the
particle -- vacancy symmetry is no longer applicable in our case.
This explains the need to develop a fundamentally new techniques
for the analysis of DDS in continuum. Despite this new techniques
cannot be applied directly in the lattice case, the embedding of
lattice systems to the continuum setting allows to obtain
(indirectly) new results for the lattice systems as well.

The paper is organized as follows. In Section~\ref{s:metric} we
introduce main statistical quantities under study: particle
densities, average velocities, etc. and derive their basic
properties. Section~\ref{s:pairing} is dedicated to the main
technical tool -- dynamical coupling. In Section~\ref{s:weak} we
apply this coupling in the weak normalization setting to prove the
uniqueness of the average velocity
(Theorem~\ref{t:velocity-density}) and to derive the complete
Fundamental Diagram for the deterministic case
(Theorem~\ref{t:fund-d-weak}). We calculate also the topological
entropy of this process (Theorem~\ref{t:entropy-exclusion}). The
strong normalization case is considered in Section~\ref{s:strong}
(Theorem~\ref{t:fund-d-strong}), while a more general setting with
local velocities of both signs is studied in
Section~\ref{s:vel-2signs}
(Theorem~\ref{t:velocity-density-2signs}). Finally, in
Section~\ref{s:gen} we discuss some generalizations of our results
and applications to certain specific traffic models.

\smallskip

{\bf Acknowledgements.} The author is grateful to B.~Gurevich,
S.~Pirogov and an anonymous referee for a number of valuable
remarks. This research has been partially supported by Russian
Foundation for Fundamental Research, Program ONIT, and French
Ministry of Education grants.

\section{Basic properties of DDS}\label{s:metric}

Here we shall study questions related to densities and velocities
of DDS. To simplify notation we use the convention that the
normalization $\cN\in\{\cNs,\cNw\}$ and specify it only if this is
necessary.

By the {\em density} $\den(x,I)$ of a configuration $x\in X$ in a
bounded segment $I=[a,b]\in\IR$ we mean the number of particles from
$x$ whose centers $x_i$ belong to $I$ divided by the Lebesgue
measure $|I|>0$ of the segment $I$. If for any sequence of {\em
nested} bounded segments $\{I_n\}$ with
$|I_n|\toas{n\to\infty}\infty$ the limit
$$\den(x):=\lim\limits_{n\to\infty}\den(x,I_n)$$
exists and does not depend on  $\{I_n\}$ we call it the {\em
density}\footnote{In Section~\ref{s:one-sided} we shall show that
   this definition may be significantly weaken in the case when
   all particles move in the same direction.} %
of the configuration $x\in X$. Otherwise one considers upper and
lower particle densities $\den_\pm(x)$ corresponding to upper and
lower limits.

The correspondence between particle densities for configurations
with $r=0$ and $r>0$ is given by the following statement.

\begin{lemma}\label{l:den-r} Let configurations $x(r)\in X(r),~r>0$ and
$x\in X$ have the same sequence of gaps $\{\Delta_i\}$. Then
$\den_\pm(x(r))=\frac{\den_\pm(x)}{1+2r\den_\pm(x)}$.
\end{lemma}

\proof Due to the one-to-one correspondence (\ref{e:r->r'}) between the
configurations $x(r)$ and $x$, for each segment $I\subset\IR^1$ which
contains $\den(x,I)\cdot|I|$ particles from the configuration $x$, one
constructs the segment $I(r)$ containing the same particles from the
configuration $x(r)$. The length of this segment is equal to
$|I(r)|=|I|+2r\cdot\den(x,I)\cdot|I|$. Therefore
$$ \den(x(r),I(r)) = \frac{\den(x,I)\cdot|I|}{|I|+2r\den(x,I)\cdot|I|}
                   = \frac{\den(x,I)}{1+2r\den(x,I)} .$$
Passing to the limit as $|I|\to\infty$ one gets the result. \qed

\begin{remark}\label{r:density} If $\exists \den(x)<\infty$ then
$|x_n-x_m|/|n-m|\toas{|n-m|\to\infty}\den(x)$.
\end{remark}

\begin{lemma}\label{l:density-preservation}
The upper/lower densities $\den_\pm(x^t)$ are preserved by
dynamics, i.e. $\den_\pm(x^t)=\den_\pm(x^{t+1})~~\forall
t\in\IZ_0$.
\end{lemma}%
%  %
\proof For a given segment $I\in\IR$ the number of particles from
the configuration $x^t\in X$ which can leave it during the next
time step cannot exceed 1 and the number of particles which can
enter this segment also cannot exceed 1. Thus the total change of
the number of particles in $I$ cannot exceed 1, because if a
particle leaves the segment through one of its ends no other
particle can enter through this end. Therefore $$|\den(x^t,I) -
\den(x^{t+1},I)|\cdot|I|\le1$$ which implies the claim. \qed%

By the (average) {\em velocity} of the $i$-th particle in the
configuration $x\in X$ at time $t>0$ we mean
$$V(x,i,t):=\frac1t\sum\limits_{s=0}^{t-1}\cN(v_i^s,x^s) %
         \equiv(x_i^{t}-x_i^0)/t.$$ %
If the limit %
$$V(x,i):=\lim\limits_{t\to\infty}V(x,i,t)$$
exists we call it the (average) {\em velocity} of the $i$-th
particle. Otherwise one considers upper and lower particle
velocities $V_\pm(x,i)$.

The correspondence between average particle velocities for configurations
with $r=0$ and $r>0$ is even simpler than for densities.

\begin{lemma}\label{l:V-r} Let configurations $x(r)\in X(r),~r>0$ and
$x\in X$ have the same sequence of gaps $\{\Delta_i\}$. Then
$\forall i,t~~V(x(r),i,t)=V(x,i,t)$ for a given collection of
local velocities $\{v_i^t\}_{i,t}$.
\end{lemma}

%  %
\proof Observe that the motion of particles depends only on the
local velocities and the sequence of gaps. Thus at any time $t\ge0$
the sequence of gaps being changing in time is still the same for both
configurations $x(r)$ and $x$. Therefore
$$ \cN(v_i^t,x^t(r)) \equiv \cN(v_i^t,x^t) ~~\forall i,t $$
which yields the claim. \qed

%Observe that the velocity may be represented as an ergodic average: %
%$V(x,i,t)=\frac1t\sum_{k=0}^{t-1}\cN(v_i^k,x^k)$.

%

\begin{lemma}\label{l:velocity-preservation}
Let $x\in X$ then $|V(x,j,t) - V(x,i,t)|\toas{t\to\infty}0$ a.s.
$\forall i,j\in\IZ$.
\end{lemma}%
\proof It is enough to prove this result for $j=i+1$. Consider the
difference between (average) velocities of consecutive particles%
\bea{ V(x,i+1,t) - V(x,i,t)
 \a= \frac{x_{i+1}^t-x_{i+1}^0}t - \frac{x_i^t-x_i^0}t \\ %
 \a= \frac{x_{i+1}^t-x_i^t}t - \frac{x_{i+1}^0-x_i^0}t \\ %
 \a= \G_i^t/t - \G_i^0/t .} %
The last term vanishes as $t\to\infty$ and it is enough to show
that the same happens with $\G_i^t/t$.

Consider first the deterministic setting (i.e. $v_i^t\equiv
v_0^t$) and show that\footnote{If $v_0^t$ takes
    both positive and negative values then %
    $\G_i^t\le\max(4v,\G_i^0)$.} %
$\forall i,t$ %
\beq{e:max-gap}{\G_i^t\le
        \function{\max(v,\G_i^0)  &\mbox{if } \cN=\cNw \\
                  \max(2v,\G_i^0) &\mbox{if } \cN=\cNs } .}%

Obviously this is true for $t=0$. Assume that this inequality
holds up to time $t\in\IZ_0$ and consider the moment $t+1$.
There might be two two possibilities: %
\begin{itemize}
\item[(a)] $\G_i^t\ge v_0^t$. Then $\cN(v_0^t,x^t)=v_0^t$ and %
$$\G_i^{t+1} = \G_i^t - \cN(v_i^t,x^t) + \cN(v_{i+1}^t,x^t)
            \le \G_i^t - v_0^t + v_{0}^t = \G_i^t
            \le \max(v,\G_i^0)$$ by the assumption. %

\item[(b)] $\G_i^t<v_0^t$. Then $\cN_w(v_0^t,x^t)=\G_i^t$ and
$\cN_s(v_0^t,x^t)=0$. Therefore %
\bea{\G_i^{t+1} \a= \G_i^t - \G_i^t + \cN(v_{i+1}^t,x^t)
            \le v \le \max(v,\G_i^0)\quad{\rm if}~\cN=\cNw ,\\
     \G_i^{t+1} \a= \G_i^t - 0 + \cN(v_{i+1}^t,x^t) \le 2v\quad{\rm
if}~\cN=\cNs.}%
\end{itemize}

Thus in the deterministic setting the gaps are uniformly bounded
in time and hence $\G_i^t/t\toas{t\to\infty}0$.

Analysis of the random setting is much more involved since the
gaps between particles in principle may grow with time and become
arbitrary large but this may happen only very slowly. To estimate
from above the value of the $i$-th gap $\G_i^t$ we drop from the
consideration all particles except the $i$-th and $(i+1)$-th
(preserving for all $t\in\IZ_0$ the velocities
$\{v_i^t,v_{i+1}^t\}_t$) and denote the resulting configuration by
$\t{x}^t:=\{\t{x}_{i}^t,\t{x}_{i+1}^t\}$ and the gap %
between this pair of particles by $\t{\G}_i^t$. We have  %
\bea{\G_i^{t+1}\a:=\G_i^t - \cN(v_i^t,x^t) + \cN(v_{i+1}^t,x^t),\\
     \t\G_i^{t+1}\a:=\t\G_i^t - \cN(v_i^t,\t{x}^t) +
 \cN(v_{i+1}^t,\t{x}^t) = \t\G_i^t - \cN(v_i^t,\t{x}^t) +
 v_{i+1}^t. } %
The comparison between $\G_i^t$ and $\t{\G}_i^t$ will be done by
induction separately for the weak and strong normalizations.

First let us prove that $\t{\G}_i^t\ge\G_i^t$ if $\cN=\cNw$. At time
$t=0$ obviously $\t{\G}_i^0=\G_i^0$. Assume that
$\t{\G}_i^t\ge\G_i^t$ for some $t\in\IZ_+$. Clearly,
$$0\le\cN(v_{i+1}^t,x^t)\le v_{i+1}^t.$$ For $v_i^t$ there
might be two possibilities: %
\begin{itemize}
\item[(a)] $v_i^t\le \G_i^t$. Then
    $\cN(v_{i}^t,x^t)=\cN(v_{i}^t,\t{x}^t)=v_i^t$ and hence %
    $$\t{\G}_i^{t+1}=\t{\G}_i^t-v_{i}^t+v_{i+1}^t
     \ge \G^t-v_{i}^t+v_{i+1}^t = \G^{t+1}.$$ %
\item[(b)] $v_i^t > \G_i^t$. Then
    $\cNw(v_{i}^t,x^t)=\G_i^t$, $\cNw(v_{i}^t,\t{x}^t)\ge\G_i^t$
    and hence %
    $$\t{\G}_i^{t+1} = \t{\G}_i^t - \cNw(v_{i}^t,\t{x}^t) + v_{i+1}^t
     \ge v_{i+1}^t = \G_i^{t+1}.$$ %
\end{itemize}

If $\cN=\cNs$ a weaker estimate $\t{\G}_i^t+v\ge\G_i^t$ takes place.
Considering again the same possibilities we see that the cases $t=0$
and (a) hold without any changes, but the case (b) should be
rewritten.

\begin{itemize}
\item[(b')] $v_i^t > \G_i^t$. Then $\cNs(v_{i}^t,x^t)=0$,
    $\cNs(v_{i}^t,\t{x}^t) =
     \function{0          &\mbox{if } v_{i}^t > \t{\G}_i^t \\
               \t{\G}_i^t &\mbox{if } v_{i}^t \le \t{\G}_i^t}$,
    and hence $\cNs(v_{i}^t,\t{x}^t)\ge\cNs(v_{i}^t,x^t)$. Thus %
    \bea{\t{\G}_i^{t+1} \a= \t{\G}_i^t - \cNs(v_{i}^t,\t{x}^t) + v_{i+1}^t \\
        \a\ge \G_i^t - v - \cNs(v_{i}^t,x^t) + v_{i+1}^t
                   - (\cNs(v_{i}^t,\t{x}^t)-\cNs(v_{i}^t,\t{x}^t))
        \ge \G_i^{t+1} - v.}
\end{itemize}

Consider now the behavior of $\t{\G}_i^t$ as a function of time
$t$. If $\t{\G}_i^t\ge v$ we get $v_i^t\le \t{\G}_i^t$ and hence
$\cN(v_{i}^t,x^t)=v_{i}^t$, which implies that outside of the
region $[0,v]$ the sequence $\t{\G}_i^t$ behave as a spatially
homogeneous reflected at $0$ random walk with i.i.d. symmetric
increments $v_{i+1}^t-v_i^t$. Thus the mathematical expectation
$E(\t{\G}_i^t)$ cannot exceed\footnote{$4v$ if local velocities
     take both positive and negative values.} %
$2v$ and hence by Chebyshev inequality the probability %
$$ P(\t{\G}_i^t/t\ge\ep) \le \frac1\ep~E(\t{\G}_i^t/t) %
                         \le \frac{2v}{t\ep} \toas{t\to\infty}0 ,$$
which finishes the proof. \qed

\?{\footnotesize On the other hand, the mathematical expectation
of the time $\t{\G}_i^t$ spends in the region $[0,v]$ is finite.
Therefore applying the Reflection Principle by the Law of Large
Numbers we obtain $\t{\G}_i^t/t\sim\sqrt{t}/t\toas{t\to\infty}0$
a.s.}

\begin{corollary}\label{c:vel-equiv}
The upper and lower particle velocities $V_\pm(x,i)$ do not depend
on $i$ (but might be random). \end{corollary}

\section{Coupling}\label{s:pairing}

Recall that a coupling of two Markov chains $x^t$ and $y^t$ acting
on the space $X$ is an arrangement of a pair of processes on
a common probability space to facilitate their direct comparison,
namely this is a pairs process $(x^t,y^t)$ defined on the direct
product space $X\times X$ satisfying the assumptions
$$P((x^t,y^t)\in A\times X)=P(x^t\in A) \quad{\rm and}\quad
  P((x^t,y^t)\in X\times A)=P(y^t\in A) $$ %
for any measurable subset $A\subseteq X$, i.e. the projections
behave as the individual processes.

Let $x^t, \2x^t$ be two copies of Markov chains, describing the
DDS which we consider throughout the paper. Typically in
continuous time interacting lattice particle systems one uses (see
e.g. \cite{Lig}) an {\em equal} coupling (pairing) when particles
sharing the same sites in the copies $x^t, \2x^t$ are considered
to be paired and all choices of their velocities are assumed to be
identical. This sort of coupling works rather well for continuous
time systems when only a single particle may move at a given
moment of time. In the discrete time case the situation is much
more complicated since an arbitrary number of particles may move
simultaneously and thus it is possible that the particles of the
processes $x^t, \2x^t$ pass each other and never share the same
positions. In fact, this difficulty is not really crucial and can
be cured under some simple technical assumptions. A more important
obstacle is that if a pair is created and only one of its members
is blocked by an unpaired particle, then due to the simultaneous
motion of the blocking unpaired particle and the non-blocked
particle belonging to the pair the following situation may happen: %
~$_{\bullet}^{\bullet\circ}\longrightarrow
  _{~~~\circ}^{~\circ~~~\circ}$. %
Thus the old pair will be destroyed but no new pair will be
created under the equal pairing construction. Here and in the
sequel we use a diagrammatic representation for coupled
configurations, where paired particles are denoted by black
circles and unpaired ones by open circles, and use the upper line
of the diagram for the $x$-particles (i.e. particles from the
$x$-process) and the lower line for the $\2x$-particles.

To deal with this obstacle we introduce a {\em dynamical}
\footnote{The word ``dynamical'' is meant to emphasize that the
   mutual arrangement of particles in pairs may change with time
   under dynamics in distinction to the conventual equal coupling
   (where the particles have coinciding positions).} %
coupling, a very preliminary version of which was described in
\cite{BP} for the lattice case and was inspired by the idea proposed
by L.~Gray for the simplest discrete time lattice TASEP
(unpublished). It is worth mention also the coupling proposed for
the lattice continuous time case by O.~Angel (see \cite{An,EFM}). As
we shall show an important advantage of the dynamical coupling with
respect to the Angel's construction is that the former guarantees
that the distances between mutually paired particles are uniformly
bounded.\footnote{In the Angel's construction the distances may grow
  to infinity.}

By the {\em dynamical coupling} of the processes $x^t,\2x^t$ we
mean a gradual pairing of close enough particles belonging to the
opposite processes satisfying the following assumptions:

\begin{itemize}

\item [(A1)] At $t=0$ all particles are assumed to be unpaired.
Velocities of mutually paired particles are identical.

\item [(A2)] Once being created a pair of particles remains
present\footnote{Starting from the moment when a pair is created
   we consider it as an entity independently on the possible
   change of particles forming it.} %
for any moment of time in the future, however at different moments
of time the roles of the pair's members may be played by different
particles.

\item [(A3)] A particle overtaking during one time step of the
dynamics some unpaired particles from the opposite process becomes
paired with one of them.

\end{itemize}

%\beq{p:overtake}{\begin{array}{l} %
%\mbox{A particle overtaking during one time step of
%the dynamics some unpaired} \\ %
%\mbox{particles from the opposite process becomes paired
%        with one of them.}
%\end{array}  } %

According to (A1)--(A3) particles from the same pair move
synchronously until either the admissibility condition breaks down
for only one of the particles (which means that its movement is
blocked by another particle) or one of the members of the pair is
swapped with an unpaired particle from the same process (see
Fig.~\ref{f:pairing} for the case of the weak normalization). It
is convenient to think about the coupled process as a ``gas'' of
single (unpaired) particles and ``dumbbells'' (pairs). A
previously paired particle may inherit the role of the unpaired
one from one of its neighbors. In order to keep track of positions
of unpaired particles we shall refer to them as $x$- and
$\2x$-{\em defects} depending on the process they belong.

%%%%%%%%%%%%%%%%%%%%%%%%%%%%%%%%%%%%%%%%%%%%%%%%%%%%%%%%%%%%%%%%%%
%% Pairing
\Bfig(200,110)
      {\put(30,95){\circle*{5}} \put(30,100){$i$}
       \bline(30,95)(1,-2)(12)
       \put(42,70){\circle*{5}}  \put(42,75){$j$}
       \put(42,65){\vector(1,0){41}} \put(44,55){$v_j=v_i$}
       \put(60,95){\circle{5}} \put(55,100){$i+1$}
       \put(60,90){\vector(1,0){10}} \put(68,80){$v_{i+1}$}
       \put(120,95){\circle{5}} \put(115,100){$i+2$}
       \put(120,90){\vector(1,0){41}} \put(135,80){$v_{i+2}$}
       \put(168,70){\circle{5}} \put(163,75){$j+1$}
       \put(167,65){\vector(1,0){8}} \put(165,55){$v_{j+1}$}
       \thicklines \bline(0,45)(1,0)(200) \thinlines
                   %\put(170,105){$t$} \put(170,30){$t+1$}
       \bezier{30}(0,95)(100,95)(200,95) \put(1,100){$x^t$}
       \bezier{30}(0,70)(100,70)(200,70) \put(1,75){$\2x^t$}
       \bezier{30}(0,20)(100,20)(200,20)    \put(1,25){$x^{t+1}$}
       \bezier{30}(0,-5)(100,-5)(200,-5)    \put(1,0){$\2x^{t+1}$}
       \put(58,20){\circle{5}} \put(57,25){$i$}
       \put(70,20){\circle*{5}} \put(70,25){$i+1$}
       \put(83,-5){\circle*{5}}  \put(83,0){$j$}
                  \bline(70,20)(1,-2)(13)
       \put(161,20){\circle*{5}} \put(155,25){$i+2$}
       \put(174,-5){\circle*{5}}  \put(173,0){$j+1$}
           \bline(161,20)(1,-2)(11) %\bline(152,25)(0,1)(100)
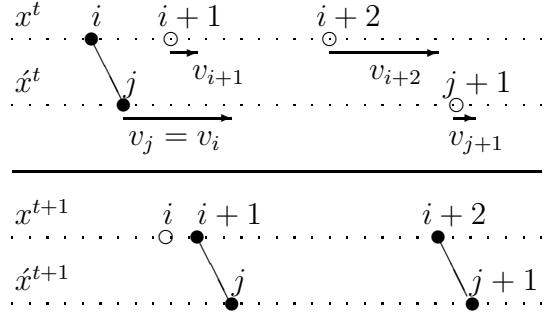
     }{Pairing of particles. Black circles corresponds to paired
        particles and open circles to defects.
        The paired particles are connected by straight lines.
        At time $t$ the particles $i$ and $j$ are paired, while
        at time $t+1$ the $x$-particle $i$ becomes unpaired and
        the $\2x$-particle $j$ becomes paired with the $x$-particle $i+1$.
        The unpaired initially particles $i+2$ and $j+1$ become
        paired at time $t+1$.
        \label{f:pairing}}
%%%%%%%%%%%%%%%%%%%%%%%%%%%%%%%%%%%%%%%%%%%%%%%%%%%%%%%%%%%%%%%%%%

There are a number of ways to realize the dynamical coupling (in
particular, using only the idea of the particle's overtaking). To
demonstrate the flexibility of our approach we describe a
different construction. Note that in the sequel we shall use only
the properties (A1)--(A3) and the proofs will not depend on other
details of the coupling.

By the $x$-{\em triple} (~$_\bullet^{~\circ~\bullet}$~ or
~$^{\bullet~\circ}_{~~~\bullet}$~) in the coupled process
$(x^t,\2x^t)$ we mean two mutually paired particles and a
$x$-defect located in the segment between them, whose index
differs by one from the index of the paired $x$-particle.
The $\2x$-triple (~$^\bullet_{~\circ~\bullet}$~ or
~$_{\bullet~\circ}^{~~~\bullet}$) is defined similarly.

Two pairs of particles are said to {\em cross} each other if
straight lines connecting positions of particles belonging to the
same pair intersect, e.g. ~$_{~\star~~\bullet}^{\bullet~~\star}$~,
where particles belonging to the same pair are marked similarly.

A $x$-defect at $x_i^t$ together with the closest\footnote{If
   there are several closest $\2x$-defects one chooses the defect
   with the smallest index.} %
$\2x$-defect at $\2x_j^t$ (~$_\circ^{~\circ}$~ or
~$^\circ_{~\circ}$~) are said to be a {\em d-pair}\/ if
$|x_i^t-\2x_j^t|<v$, this pair of defects does not cross with any
mutually paired particles, and the open segment $(x_i^t,\2x_j^t)$
does not contain any other defects. We say that a d-pair $(i,j)$
is {\em smaller} than a d-pair $(n,m)$ if $|i|<|n|$, or if $i<n$
in case $|i|=|n|$. Observe that $i=n$ but $j\ne{m}$ cannot happen
in distinction to $i\ne{n}$ but $j=m$.

Note that in the collection
~$_{\bullet~\bullet}^{~\circ~\bullet~\bullet}$~ the first two
$x$-particles together with the first $\2x$-particle form a
$x$-triple despite the presence of an additional paired particle
in the segment between them. On the other hand, the collection
~$_{\circ~\bullet}^{~\bullet~\circ}$~ does not contain neither
triples nor d-pairs.

A pair of configurations $(x^t,\2x^t)$ representing the coupled
process at time $t$ is said to be {\em proper} if it does not
contain $x$- or $\2x$-triples, d-pairs, and crossing mutually
paired particles.

The fact that at time $t$ the pair of configurations $(x^t,\2x^t)$
were proper does not imply that it remains proper under dynamics
at time $(t+1)$. In particular, triples of both types and d-pairs
may be created, e.g.
~$_{\bullet~~~\circ}^{~~~\bullet}\longrightarrow
  _{\bullet~\circ}^{~~~\bullet}$~ or %
$~_{~~~~\circ\circ}^\circ~\longrightarrow~_{\circ\circ}^\circ$,
however due to the particle order preservation crossing mutually
paired particles cannot appear.

\begin{lemma}\label{l:triples} Let a pair of configurations
$(x^t,\2x^t)$ have no crossing mutually paired particles. Then
among triples of the same kind there are no common elements.
\end{lemma}
\proof Direct inspection. As an illustration let us check the claim
about $x$-triples. Assume that two $x$-triples have a common
$x$-defect (mutually paired particles cannot be common by definition).
Then this implies that the mutually paired particles in these
triples either cross each other
$~_{~~\star~~\bullet}^{\bullet~~\circ~~\star}~$
or the index of one of the paired $x$-particles differs from the
index of the common defect by more than one
$~_{~~~~~\star~\bullet}^{\star~\bullet~\circ}~$. The latter
contradicts to the definition of the $x$-triple, why the former
contradicts to the assumption about the absence of crossing
mutually paired particles. In the diagrams above paired particles
from the 2nd triple are marked by stars to distinguish them from
the 1st triple. \qed

Therefore all triples of the same kind may be resolved
simultaneously.
This will be done as follows. A $x$- or $\2x$-triple is
transformed such that the former defect is becoming paired to the
particle from another process, while another previously paired
particle is becoming unpaired:
~$_\bullet^{~\circ~\bullet}\longrightarrow~_\bullet^{~\bullet~\circ}$~.

The case of a d-pair is even simpler, namely the defects
``annihilate'' forming mutually paired particles:
~$_\circ^{~\circ}\longrightarrow~_\bullet^{~\bullet}$~. In all
cases the positions of particles are preserved but their ``roles''
are changing.

Finally the coupling procedure consists of the following steps: %
\begin{itemize}
\item[(1)] Each $x$-triple is recursively resolved:
  ~$_\bullet^{~\circ~\bullet}\longrightarrow~_\bullet^{~\bullet~\circ}$~.
\item[(2)] Each $\2x$-triple is recursively resolved:
  ~$^\bullet_{~\circ~\bullet}\longrightarrow~^\bullet_{~\bullet~\circ}$~.
\item[(3)] The smallest\footnote{The ordering of d-pairs is
  updated after each recursion procedure.} %
  d-pair is recursively resolved:
   ~$_\circ^{~\circ}\longrightarrow~_\bullet^{~\bullet}$~.
\end{itemize}

\begin{lemma}\label{l:coupling-markov} The coupling procedure
described above is well defined, leads to the Markovian coupling,
and satisfies the assumptions (A1)--(A3).
\end{lemma}
\proof Let us check that this procedure is well defined. By
Lemma~\ref{l:triples} if a particle belongs to a certain triple
then it cannot belong to any other triple. On the other hand,
segments belonging to paired particles may overlap and resolving a
$x$- or $\2x$-triple one may create a new one of the same kind:
$$_{\bullet~~\bullet}^{~\circ~~\bullet~\bullet}\longrightarrow
 ~_{\bullet~~\bullet}^{~\bullet~~\circ~\bullet}\longrightarrow
 ~_{\bullet~~\bullet}^{~\bullet~~\bullet~\circ}.$$ %
This explains the necessity of the recursion during the first two
steps of the procedure. Note that resolving a $x$-triple one
cannot create a new $\2x$-triple and vice versa (defects do not
move from one process to another).

Elements of the smallest d-pair might belong to some other
d-pairs. Therefore resolving it we might change the d-order of the
remaining d-pairs. To take this into account we are recalculating
the d-order after each recursion procedure.

Consider now the motion of a given defect under the recursions in
the coupling procedure. Observe that the defect may move arbitrary
far in any direction from its initial position due
to these recursions: %
$$_{~\circ~\bullet~\bullet~\cdots~\bullet~\bullet}
     ^{\bullet~\bullet~\cdots~\bullet~\bullet}\longrightarrow~
     _{~\bullet~\bullet~\cdots~\bullet~\bullet~\circ}
     ^{\bullet~\bullet~\cdots~\bullet~\bullet}.$$
Nevertheless a defect cannot change its direction of
movement. Assume from the contrary that a $x$-defect during two
consequent steps of the recursion moved first to the right
($_\bullet^{~\circ~\bullet}\longrightarrow~_\bullet^{~\bullet~\circ}$)
and then to the left
($_{~~~\bullet}^{\bullet~\circ}\longrightarrow
 ~_{~~~\bullet}^{\circ~\bullet}$). %
This can happen only if after the first step of the recursion the
defect became a member of a new $x$-triple of type
$_\bullet^{~\circ~\bullet}$. Then the only candidate for the role
of the paired $x$-particle in this $x$-triple is the paired
$x$-particle which played the role of this defect on the previous
recursion step. We came to the contradiction, because a particle
may belong to only one pair.

Thus the recursion is finite in the sense that each defect in a
bounded spatial segment in finite time either will stop moving or
will leave this segment and never return back. Note however that
in general one cannot divide a configuration into finite pieces
and deal with them separately since a defect may move from one
piece to another.

After the application of the first two steps all $x$- or
$\2x$-triples will be eliminated and only d-pairs may be present.
Observe now that when one resolves a d-pair neither triples nor
new defects are created. However since various d-pairs may
intersect they should be resolved separately during the last step.
Additionally neither of above procedures may create crossing pairs
of mutually paired particles (since members of different triples
of the same type do not intersect and c- and d-pairs cannot cross
each other).

Let the pair of configurations $(x^{t-1},\2x^{t-1})$ be proper.
Then according to arguments above after one time step of the
dynamics the application of the coupling procedure, is well
defined and the pair of configurations $(x^t,\2x^t)$ at time $t$
is proper as well.

By the construction the one-time step transition probabilities for
both processes $x^t$ and $\2x^t$ remain unchanged and the one-time
step transition probabilities for the pairs process are well
defined. Therefore this construction defines a Markovian coupling
between two copies of the Markov chain describing our DDS.

The property (A1) holds by the construction. A pair breaks down
only if one of its members is replaced by an unpaired particle,
and hence the pair as a whole survives. This proves (A2). The
property (A3) follows from the fact that under the one time step
of the dynamics of a proper pair of configurations all objects
under consideration: $x$- and $\2x$-triples, and d-pairs may be
created only during the particles overtaking. \qed

Denote by $\den_u(x,I)$ the density of the $x$-defects belonging
to a finite segment $I$, and by $\den_u(x):=\den_u(x,\IR)$ the
upper limit of $\den_u(x,I_n)$ taken over {\em all} possible
collections of nested finite segments $I_n$ whose lengths go to
infinity.

%\n{\bf Definition.} %
We say that a coupling of two Markov particle processes
$x^t,\2x^t$ is {\em nearly successful} if the upper density of the
$x$-defects $\den_u(x)$ vanishes with time a.s.
This definition differs significantly from the conventional
definition of the successful coupling (see e.e. \cite{Lig}), which
basically means that the coupled processes converge to each other
in finite time.

In the random setting under some regularity assumptions the
dynamical coupling turns out to be nearly successful (the proof of
this result goes out of the scope of the present paper and will be
published elsewhere), however in general especially in the
deterministic setting this property needs not hold.

Applying the notion of the nearly successful coupling to the
exclusion process under study we get the following conditional
result.

\begin{lemma}\label{l:velocity-coupled} Let $x,\2x\in X$ with
$\den(x)=\den(\2x)$, and let there exist a nearly successful
coupling $(x^t,\2x^t)$ such that distances between the pair members
are uniformly bounded from above by $\gamma(t)=o(t)$. Then
$$|V(x,0,t)-V(\2x,0,t)|\toas{t\to\infty}0.$$
\end{lemma}%
\proof Consider an integer valued function $n_t$ which is equal to
the index of the $\2x$-particle paired at time $t>0$ with the
$0$-th $x$-particle. If the $0$-th $x$-particle is not paired at
time $t$ we set %
$n_t:=\function{n_{t-1} &\mbox{if } t>0 \\
                0       &\mbox{if } t=0}$. %

To estimate the growth rate of $|n_t|$ at large $t$ observe that
$n_t$ changes its value only at those moments of time when the
$0$-th $x$-particle meets a $\2x$-defect. By the assumption about
the nearly successful coupling at time $t\gg1$ the average distance
between the defects at time $t$ is of order $1/\den_u(\2x^{t})$
while the amount of time needed for two particles separated by the
distance $L$ to meet cannot be smaller than $L/(2v)$. Therefore the
frequency of interactions of the $0$-th $x$-particle with
$\2x$-defects may be estimated from above by the quantity of order %
$\den_u(\2x^t)\toas{t\to\infty}0$, which implies %
$n_t/t\toas{t\to\infty}0$.

Now we are ready to prove the main claim. %
\bea{ |V(x,0,t) - V(\2x,0,t)| \a= |(x_0^t-x_0^0) -
(\2x_0^t-\2x_0^0)|/t \\
\a\le |x_0^t-\2x_0^t|/t + |x_0^0-\2x_0^0|/t \\%
\a\le |x_0^t-\2x_{n_t}^t|/t +
        \frac{|n_t|}t~|\2x_{n_t}^t-\2x_0^t|/|n_t|
                            + |x_0^0-\2x_0^0|/t .} %
The 1st addend can be estimated from above by
$\gamma(t)/t\toas{t\to\infty}0$. The 2nd addend is a product of
two terms $|n_t|/t$ and $|\2x_{n_t}^t-\2x_0^t|/|n_t|$. As we have
shown, the 1st of them vanishes with time. If $|n_t|$ is uniformly
bounded, then the 2nd term is obviously uniformly bounded on $t$.
Otherwise, for large $|n_t|$ by Remark~\ref{r:density} and the
density preservation the 2nd term is of order $\rho(\2x)$, which
proves its uniform boundedness as well. Thus the 2nd addend goes
to 0 as $t\to\infty$. Noting finally that the last addend also
vanishes with time we are getting the result. \qed

\section{Weak normalization}\label{s:weak}

Consider the coupled process $(x^t,\2x^t)$ under the weak
normalization and set $W_{ij}^t:=x_i^t-\2x_j^t$.

\begin{lemma}\label{l:dist-pair}
The supremum of ~$|W_{ij}^t|$~ taken over all mutually paired
particles is uniformly bounded by $v$ for any $t\in\IZ_0$.
\end{lemma}
\proof We start at time $t=0$ when there are no pairs and wait
until the first of them appears. At that moment the distance
between the members in a pair cannot exceed $v$. Starting from
that moment the distances may grow and some new pairs may be
created. Contrary to our claim assume that there is the first
moment of time $t$ at which there is a pair of particles located
at $x_i^t,\2x_j^t$ for which $|x_i^t - \2x_j^t|>v$ and it is the
largest distance between the paired particles at that moment of
time (or one of the largest) and such that $|x_i^{t-1} -
\2x_j^{t-1}|\le v$. According to the definition of the pairing
process there are no unpaired particles between the particles from
the same pair. Therefore in order to enlarge the distance between
the particles one of them should be blocked by a particle from
another pair, which contradicts to the assumption about the
maximality of the distance. \qed

\begin{lemma}\label{l:s-coupling}
Let $\den(x)=\den(\2x)$ and let in the coupled process $\forall
i,j~~\exists$ a (random) moment of time $t_{ij}<\infty$ such that
$x_i^t>\2x_j^t$ for each $t\ge t_{ij}$. Then the coupling is nearly
successful.
\end{lemma}
\proof By the assumption each $x$-particle will overtake
eventually each $\2x$-particle located originally to the right
from its own position and thus will form a pair with it or with
one of its neighbors (if they are so close that were overtaken
simultaneously). Thus the creation of pairs is unavoidable. To
show that the upper density of defects cannot remain positive,
consider how the defects move under our assumptions. Assume that
at time $t\ge0$ the $i$-th $x$-particle is paired with the $j$-th
$\2x$-particle. Then by Lemma~\ref{l:dist-pair} in order to
overtake at time $s>t$ the $j$-th $\2x$-particle significantly (by
a distance larger than $v$) the $i$-th $x$-particle necessarily
needs to break the pairing with the $j$-th $\2x$-particle. Thus by
the property (A3) of the dynamical coupling either a $x$-defect
overtakes the
$j$-th $\2x$-particle: %
~$_{~\bullet}^{\circ~~~\bullet}~\longrightarrow~
  _\bullet^{~\circ~\bullet}~\longrightarrow~
  _\circ^{~\bullet~\circ}$,
or the $i$-th $x$-particle overtakes a $\2x$-defect:
~$_{\bullet~~~\circ}^{~~~\bullet}\longrightarrow~
  _{\bullet~\circ}^{~~~\bullet}\longrightarrow~
  _{\circ~\bullet}^{~~~\bullet}$. %
(Otherwise this pair will not be broken.) Therefore during this
process the $x$-defects move to the right while the $\2x$-defects
move to the left. Hence they inevitably meet each other and
``annihilate''. The assumption about the equality of particle
densities implies the result. \qed

\subsection{Uniqueness of the average velocity}
\label{s:vel-uniq}

As we shall see under our assumptions even in the weak
normalization case the nearly successful coupling needs not hold
(e.g. in the deterministic setting). Therefore one cannot apply
directly Lemma~\ref{l:velocity-coupled} in this case. Nevertheless
we shall show that the absence of coupling is not a serious
obstacle and it can be used as a diagnostic tool.

\begin{theorem}\label{t:velocity-density} %
In the weak normalization case the set of limit points as
$t\to\infty$ of the sequence $\{V(x,t)\}_{t\in\IZ_0}$ depends only
on the density $\den(x)$ assuming that the latter is well defined.
\end{theorem}%
\proof Consider a general DDS under the weak normalization. Let
$x,\2x\in X_\den:=\{z\in X:~~\den(z)=\den\}$ be two admissible
configurations of the same particle density. If one assumes that
the coupling procedure described in Section~\ref{s:pairing} leads
to the nearly successful coupling of particles in these
configurations then by Lemma~\ref{l:dist-pair} the assumptions of
Lemma~\ref{l:velocity-coupled} are satisfied and hence
$|V(x,0,t)-V(\2x,0,t)|\toas{t\to\infty}0$ which by
Lemma~\ref{l:velocity-preservation} implies the claim.
In general the assumption about the nearly successful coupling may
not hold,\footnote{Consider e.g. the deterministic
     setting with $1/\den>5v$ and the configurations $x_i:=i/\den$
     and $\2x_i:=i/\den+2v$. Then $\den(x)=\den(\2x)=\den$,
     $V(x)=V(\2x)=v$ but no pair will be created.} %
however as we demonstrate below the pairing construction is still
applicable.

\?{The idea is based on the observation that if one assumes that
the claim does not hold then we are in a position to apply
Lemma~\ref{l:s-coupling} which implies the nearly successful
coupling which by Lemma~\ref{l:velocity-coupled} contradicts to
the assumption.}

Define random variables
$$W_{ij}^t:=x_i^t-\2x_j^t,~i,j\in\IZ,~t\in\IZ_0.$$ Then %
$$V(x,i,t)-V(\2x,j,t)=W_{ij}^t/t - W_{ij}^0/t.$$ Since by
Lemma~\ref{l:velocity-preservation} the differences between
average velocities of different particles belonging to the same
configuration vanish with time it is enough to consider only the
case $i=j=0$. For $W_{00}^t$ there might be three possibilities
which we study separately:

\begin{itemize}
\item[(a)] $\lim\limits_{t\to\infty} W_{00}^t/t=0$. Then %
$|V(x,0,t)-V(\2x,0,t)|\le|W_{00}^t|/t + |W_{00}^0|/t
 \toas{t\to\infty}0$, %
which by Corollary~\ref{c:vel-equiv} implies that the sets of
limit points of the average velocities coincide.

\item[(b)] $\limsup\limits_{t\to\infty}W_{00}^t/t>0$. Then
$\forall i\in\IZ$ the $i$-th particle of the $x$-process will
overtake eventually each particle of the $\2x$-process located at
time $t=0$ to the right from the point $x_i^0$. This together with
the assumption of the equality of particle densities allows to
apply Lemma~\ref{l:s-coupling} according to which the coupling is
nearly successful. On the other hand, by Lemma~\ref{l:dist-pair}
the distance between mutually paired particles cannot exceed $v$.
Therefore by
Lemma~\ref{l:velocity-coupled} we have ~~%
$|V(x,0,t)-V(\2x,0,t)|\toas{t\to\infty}0$, %
which contradicts to the assumption (b).

\item[(c)] $\limsup\limits_{t\to\infty} W_{00}^t/t<0$. Changing
the roles of the processes $x^t,\2x^t$ one reduces this case to
the case (b).
\end{itemize}

Thus only the case (a) may take place. \qed

\subsection{Deterministic setting}

\begin{theorem}\label{t:fund-d-weak} (Fundamental Diagram)
In the deterministic setting %
\beq{e:FD-w}{V(x)=\lim\limits_{t\to\infty}\frac1t
         \sum_{s=0}^{t-1}\min(1/\rho,v_0^s)
=\function{v &\mbox{if ~~ } \den(x)\le 1/v  \\
           1/\den(x) &\mbox{otherwise }} } %
if $v_0^t\equiv v$.
\end{theorem}
\proof Consider a family %
$$\0X_\rho:=\{x\in X:~~x_i:=i/\rho+\omega,~\omega\in\IR\}$$ %
of uniformly spatially distributed configurations of a given
density $\rho>0$. This set is forward invariant and
$$x_i^{t+1}-x_i^{t}\equiv\min(1/\rho,v_0^t)~~\forall
x^t\in\0X_\rho, i\in\IZ,$$ %
i.e. all particles in the configuration get the same normalized
local velocity $\min(1/\rho,v_0^t)$
(depending in general on time $t$). %
By the definition of the deterministic setting the limit %
$$V(x):=\lim\limits_{t\to\infty}\frac1t
         \sum_{s=0}^{t-1}\min(1/\rho,v_0^s)$$ %
is well defined. On the other hand, by
Theorem~\ref{t:velocity-density} all configurations of the same
density have the same average velocity, which implies the result.
\qed

\begin{remark} This result looks very similar to the one known for
the deterministic version of the lattice TASEP (see
\cite{NS,Bl-erg}), however the latter case is characterized by the
following feature: if the density is large enough particles
inevitably form dense clusters without vacancies inside (static
traffic jams). The proof above shows that the ``typical'' behavior
of high density configurations in continuum is different: they do
form particle clusters, but these clusters are not staying at rest
but are moving at a constant velocity as an ``echelon''. It is of
interest that in order to imitate such behavior a number of
complicated lattice models were developed.
\end{remark}

\begin{remark} The construction used in the proof is especially
striking in that the same family of uniformly spatially
distributed configurations allows to study the limit dynamics in
the deterministic setting for all configurations having densities.
Note that this argument cannot be applied directly in the lattice
version of DDS. Nevertheless since the ``lattice configurations''
are included in DDS under consideration the result holds as well,
which implies completely new results for lattice TASEPs with long
jumps. \end{remark}

\begin{corollary}\label{c:V-r-w} Let $x(r)\in X(r),~r>0$ and
$\den(x(r))$ be well defined and let ~$\forall i,t~v_i^t\equiv v$. Then
$$ V(x(r))=\function{v &\mbox{if ~~ } \den(x)\le \frac1{v+2r}  \\
           1/\den(x(r))-2r &\mbox{otherwise }}.$$
In particular in the lattice setting this reads
$$ V(x(1/2))=\function{v &\mbox{if ~~ } \den(x)\le \frac1{v+1}  \\
           1/\den(x(1/2))-1 &\mbox{otherwise }}.$$
\end{corollary}
\proof By (\ref{e:r->r'}) and Lemma~\ref{l:den-r} for each configuration
$x(r)$ one constructs the configuration $x$ with the same sequence of gaps
and the relation between their densities is written as
$$ \den(x)=\frac{\den(x(r))}{1-2r\den(x(r))} .$$
Additionally by Lemma~\ref{l:V-r} average velocities related to
configurations with the same sequence of gaps coincide. Substituting
$\den(x)$ as a function of $\den(x(r))$ to (\ref{e:FD-w}) we get
the result. \qed

\subsection{Entropy}\label{s:entropy}

In this Section we restrict the analysis to the pure deterministic
setting (i.e. $v_i^t\equiv v~~\forall i,t$). Then our DDS is
defined by a deterministic map $\map_v:X\to X$ from the set of
admissible configurations into itself. Our aim is to show that
this map is chaotic in the sense that its topological entropy is
infinite.\footnote{Normally one says that a map is chaotic if
   its topological entropy is positive, so infinite value of
   the entropy indicates a very high level of chaoticity.} %

We refer the reader to \cite{Bi,Wa} for detailed definitions of
the topological and metric entropies for deterministic dynamical
systems and their properties that we use here. To avoid
difficulties related to the non-compactness of the phase space we
define the topological entropy of a map $\map_v$ (notation
$h_{{\rm top}}(\map_v)$) as the supremum of metric entropies of
this map taken over all probabilistic invariant measures (compare
to the conventional definition of the topological entropy and its
properties in \cite{Wa}).

For a finite subset of integers $I$ and a collection
$C:=\{C_i\}_{i\in I}$ of open intervals the subset %
$\cC_{I,C}:=\{x\in X:~~x_i\in C_i~~\forall i\in I\}$ is called a
finite {\em cylinder}.\footnote{In general the cylinder $\cC_{I,C}$
   might be empty for nonempty sets $I,C$.} %
We endow the space of admissible configurations $X$ by the
$\sigma$-algebra $\cB$ generated by the finite cylinders defining a
topology in this space.

We start the analysis with the action of a shift-map in continuum
$\sigma_v:X\to X$ defined as %
$$(\sigma_v x)_i:=x_i+v~~~ i\in\IZ, x\in X.$$

\begin{lemma}\label{l:entropy-shift-map}
The topological entropy of the shift-map in continuum $\sigma_v$
is infinite.
\end{lemma}
\proof The preimage of a finite cylinder under the action of
$\sigma_v$ is again a finite cylinder. Therefore this map is
continuous in the topology induced by the $\sigma$-algebra $\cB$
generated by finite cylinders.

The idea of the proof is to construct an invariant subset of $X$
on which the map $\sigma_v$ is isomorphic to the full shift-map in
the space of sequences with a countable alphabet. The result
follows from the observation that the topological entropy of the
full shift-map $\sigma^{(n)}$ with the alphabet consisting of $n$
elements is equal to $\ln n$ (see, e.g. \cite{Bi,Wa}).

Let $\alpha:=\{\alpha_i\}_{i\in\IZ_+}$ with $\alpha_i\in(0,v)$ and
let $\alpha^n:=\{\alpha_i\}_{i=1}^n$. Consider a sequence of
subsets $X^{(n)}\subset X$ consisting of {\em all} configurations
$x\in X$ satisfying the condition $\forall k\in\IZ~~
x_{2k}\in{v}\IZ, ~x_{2k+1}\in x_{2k} + \alpha^n$. Then $X^{(n)}$
is $\sigma_v$-invariant and the restriction $\sigma_v|X^{(n)}$ is
isomorphic to the full shift-map $\sigma^{(n)}$ with the alphabet
$A^n$ consisting of $n$ elements $\{a_i\}$ of type %
$a_i:=\{[0,\alpha_i),[\alpha_i,v)\}$, i.e. each element is
represented by a pair of neighboring intervals. Therefore the
topological entropy of $\sigma^{(n)}$ is equal to %
$\ln n\toas{n\to\infty}\infty$. \qed

Another elegant (but technically difficult) way to derive this
result was proposed by Boris Gurevich. Consider a special flow $S^t$
corresponding to the shift-map acting on the sequences $\{\G_i(x)\}$
with the roof function equal to the first nonnegative particle
coordinate. This shift-map has an infinite alphabet, hence its
entropy is infinite. The special flow $S^1$ is isomorphic to the
1-shift of $\{x_i\}$, while the entropy of the special flow can be
calculated by the Abramov-Rohlin formula.

\begin{theorem}\label{t:entropy-exclusion}
The topological entropy of the pure deterministic exclusion
process in continuum is infinite.
\end{theorem}
\proof The preimage of a finite cylinder under the action of
$\map_v$ is again a finite cylinder. Therefore this map is
continuous in the topology induced by the $\sigma$-algebra $\cB$
generated by finite cylinders.

Observe that the subset %
$X_0:=\{x\in X:~~\G_i(x)\ge v~~ \forall i\in\IZ\}$ %
of the set of admissible configurations is $\map_v$-invariant.
Therefore $h_{{\rm top}}(\map_v)\ge h_{{\rm top}}(\map_v|X_0)$
and for our purposes it is enough to show that the latter is
infinite. On the other hand, by the definition of the map
$\map_v$ we have $\map_v|X_0\equiv\sigma_v|X_0$.

We still cannot apply the result of
Lemma~\ref{l:entropy-shift-map} directly because in the case under
consideration the gaps between particles are greater or equal to
$v$ by the construction. Recall that in the proof of
Lemma~\ref{l:entropy-shift-map} the gaps were not greater than
$v$. To this end one sets $\alpha_i\in(v,2v)$ and
modifies the definition of $X^{(n)}$ as follows: %
$$x_{2k+1}\in x_{2k} + \alpha^n\quad \forall k\in\IZ,~
              x_{2k}\in3v\IZ .$$ %
Consider the the alphabet $A^{(n)}$ with elements of type
$a_i:=\{[0,\alpha_i),[\alpha_i,3v)\}$. Then the $3$-d power of the
map $\map_v|X_0$ is isomorphic to the full shift-map
$\sigma^{(n)}$
with the alphabet  $A^{(n)}$. Using that %
$$3h_{{\rm top}}(\map_v|X_0)=h_{{\rm top}}((\map_v|X_0)^3)
=h_{{\rm top}}(\sigma^{(n)})=\ln n$$ we get the result. \qed

\section{Strong normalization}\label{s:strong} %

Recall that $W_{ij}^t:=x_i^t-\2x_j^t$ for $x^t,\2x^t\in X,~t\ge0$.

\begin{lemma}\label{l:gap-strong}
There exists a coupled process $(x^t,\2x^t)$ such that under the
strong normalization $\sup_{i,j,t}W_{ij}^t=\infty$, where the
supremum is taken over all mutually paired particles.
\end{lemma} %
\proof It seems that the argument applied in the weak
normalization case should work also in the case of the strong
normalization. However, a close look shows that in this case a
``blocked'' particle does not move to ``touch'' the particle
conflicting with it (as it would in the weak normalization case)
but preserves its position instead.
Therefore the distance between members of the same pair may become
larger than the distance between the members of the ``blocking''
pair which cannot happen in the weak
normalization case: %
~$_{\bullet~~\bullet~~~~\bullet~~~}
  ^{~~\bullet~~~~\bullet\bullet~~~}\longrightarrow~
  _{\bullet~~~~~\bullet~~~\bullet}
  ^{~~~~~\bullet\bullet~~~\bullet}$. %
Here initially distances between members in pairs do not exceed
$v$. The 1st pair is blocked by the 2nd pair and since the
$\2x$-member of the 1st pair cannot move (while the $x$-member
can) the distance between them becomes larger than $v$.

To demonstrate that distances between members in pairs may grow to
infinity fix some $0<\ep\ll1$ and consider a pair of
configurations $x,\2x$ such that $x_0=\2x_0=0$ and
$\G_{2k}=\frac32(v-\ep), ~\G_{2k+1}=\frac12(v-\ep),
~\2\G_k=v-\ep~\forall k\in\IZ$. After the application of the
pairing procedure $\forall i$ the $i$-particles in both
configurations will become paired forever. On the other hand,
under dynamics $\2x^t\equiv\2x^0~\forall t$ while the
$x$-particles having gaps greater than $v$ will at constant
velocity $v$. Therefore the distances between members in pairs
will grow linearly with time. \qed

\?{Nevertheless for spatially periodic configurations
$\sup_{ij}W_{ij}^t<\infty$ if density $\den<1/v$. Indeed the small
density implies the presence of gaps $>v$. Hence the proof of
Lemma~\ref{l:gap-strong} breaks down. On the other hand, due to
the spatial periodicity it is enough to consider only one spatial
period, for which the claim is trivial.}

This result demonstrates and partially explains a significant
difference in the behavior of DDS under weak and strong
normalizations. Still, as we are going to show, at least some
features of the Fundamental Diagram are preserved. Consider the pure
deterministic setting (i.e. $v_i^t\equiv v$).
The inequality~(\ref{e:max-gap}) shows that in this case gaps
between particles cannot become much larger than their initial
values. The following result demonstrates that under some mild
additional assumptions (which definitely hold for high particle
densities) large gaps will disappear with time.

\begin{lemma}\label{l:max-gap} Let $x\in X$ be spatially periodic
and we consider only the pure deterministic setting (i.e.
$v_i^t\equiv v$). Assume that $\forall t~\exists j>t:~\G_j(x^t)<v$.
Then $\forall i~\exists t_i<\infty:~\G_i(x^t)<2v~\forall t\ge t_i$.
\end{lemma}
\proof Observe that the spatial periodicity and its period is
preserved under the pure deterministic dynamics. Thus the situation
is equivalent to the consideration of a finite number (say $N$)
particles on a ring and to the assumption that for each $t\in\IZ_+$
among these particles there is a particle with a gap less than $v$
ahead of it. Note that according to the definition of the strong
normalization $\cN_s(v_i^t,x^t)=0$ whenever $\G_i(x^t)<v$.
By (\ref{e:max-gap}) $\G_i(x^t)<2v$ implies $\G_i(x^{t+1})<2v$.
Therefore new new long gaps (of size larger or equal to $2v$)
cannot be created and we need to show only that long gaps in the
original configuration will cease to exist with time.

By the assumption for any $t$ there exists a short gap (shorter
than $v$) and the corresponding particle will not move during the
next time step. Thus the index of the short gap decreases by one
after each time step until it ``collides'' with one of the long
gaps: $\G_i(x^t)\ge2v, ~\G_{i+1}(x^t)<v$. On the next time step
$\G_i(x^{t+1}):=\G_i(x^t)-v$. Due the spatial periodicity the
amount of time between these ``collisions'' is bounded and after
each of them the length of a long gap decreases by $v$. Thus they
will disappear in finite time. \qed

%%%%%%%%%%%%%%%%%%%%%%%%%%%%%%%%%%%%%%%%%%%%%%%%%%%%%%%%%%%%%%%%%%
%% Fundamental diagram
\Bfig(150,100)
      {%\footnotesize{
       \put(0,0){\vector(1,0){140}} \put(0,0){\vector(0,1){90}}
       \thicklines
       \bline(0,73)(1,0)(60)
       \bezier{200}(60,73)(75,5)(125,2)
       \bezier{200}(30,73)(45,10)(60,0)
       \thinlines
       \bezier{30}(30,73)(30,36)(30,0)
       \bezier{30}(60,73)(60,36)(60,0)
       \put(145,0){$\rho$}  \put(-8,85){$V$} \put(-8,70){$v$}
       \put(85,30){$\frac1\rho$} \put(17,15){$\frac1\rho-v$}
       \put(-8,-5){$0$}  \put(54,30){$H$}
       \put(25,-12){$\frac1{2v}$} \put(57,-12){$\frac1{v}$}
       %\bline(34,61)(1,1)(10) \bline(35,55)(1,1)(16)
      %}
      }
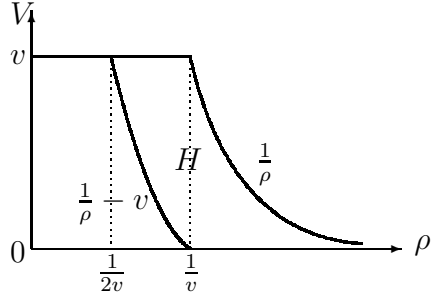
{Fundamental Diagram (dependence of the average velocity $V$ on
the particle density $\rho$) for the pure deterministic setting
under the strong normalization. The curvilinear region %
$H:=\{(\den,V): ~\frac1\rho-v\le V\le\frac1\den, ~V\le v\}$ %
corresponds to the hysteresis phase. \label{f:velocity-strong}}
%%%%%%%%%%%%%%%%%%%%%%%%%%%%%%%%%%%%%%%%%%%%%%%%%%%%%%%%%%%%%%%%%

%%

\begin{theorem}\label{t:fund-d-strong} Let $x\in X$ and $\den(x)$
be well defined. Then $V(x)=v$ if $\den(x)<\frac1{2v}$ and otherwise
for a.e. point $(\den,V)$ in the curvilinear region %\newline%
$$H:=\{(\den,V): ~\max(1/\rho-v,0)\le V\le\min(1/\den,v)\}$$ %
(see Fig.~\ref{f:velocity-strong}) %
there exists a configuration $x\in X$ with $\den(x)=\den, V(x)=V$,
i.e. the region $H$ corresponds to the hysteresis.
\end{theorem} %

\proof We say that particles numbered from $i+1$ to $i+k$ with
$i\in\IZ, k\in\IZ_+$ belonging to an admissible configuration $x\in
X$ form a {\em cluster} of {\em length} $k$ if all gaps between them
are strictly less than $v$ and the gaps to surrounding particles are
not smaller than $v$, i.e. $\G_{i+j}< v~\forall j=1,2,\dots,k-1$ and
$\G_i,\G_{i+k}\ge v$. Positions of particles belonging to the
cluster are changing with time, and leading particles leave it,
while some new particles may join the cluster from the other side.
Nevertheless the length of a cluster cannot grow with time (and new
clusters cannot be born in the pure deterministic setting in
distinction to the random one) since the rate with which the leading
particle leaves the cluster (one per unit time) is at least not
smaller than the rate at which new particles join the cluster from
the other side.

We start with the analysis of configurations of low density
(smaller than $\frac1{2v}$) and our aim is to show that in this
case each particle achieves eventually the largest available
velocity $v$.
Consider the motion of the $0$-th particle in a configuration
$x\in X$ with $0<\den(x)<\frac1{2v}$ and denote by $\hat{t}$ the
first moment of time after which this particle will not join any
cluster. If $\hat{t}<\infty$ then
$\cN_s(v_0^t)\equiv~\forall t\ge\hat{t}$ and hence %
$V(x,0,t)\toas{t\to\infty}v$.

If $\hat{t}=\infty$ then there exists an infinite sequence of
clusters of growing length such that the $0$-th particle joins
each of them consecutively. Let us show that this assumption
contradicts to the condition that $\den(x)<\frac1{2v}$. We number
the clusters to which the $0$-th particle will join according to
their natural order starting from $k=1$ and introduce the
following notation: $t_k$ -- the moment of time when the $0$-th
particle joins the $k$-th cluster, $n_k$ -- the number of
particles in this cluster, $m_k$ -- the number of particles in the
open segment between $x_0$ and the beginning of this cluster, and
$L_k$ -- the length of the minimal segment containing the $k$-th
cluster and the point $x_0$. Then
$$ \den(x,(x_0,x_0+L_k]) =
   \frac{m_k+n_k}{L_k}\toas{k\to\infty}\den(x) .$$
All $m_k$ particles will join the $k$-th cluster during the time
$t_k$ and at time $t_k$ this cluster should still exist. Therefore
the distance which the $0$-th particle covers during this time
cannot be smaller than $L_k-m_kv-n_kv$ while its velocity cannot
exceed $v$ and thus %
$$ t_kv \ge L_k-m_kv-n_kv .$$
On the other hand, exactly $t_k$ particles will leave the cluster
during this time, i.e. $m_k+n_k\ge t_k$. This gives %
\beq{e:den}{ \frac{m_k+n_k}{L_k} \ge \frac{t_k}{L_k}
  \ge \frac{L_k/v - m_k - n_k}{L_k}
    = \frac1v - \frac{m_k+n_k}{L_k} .} %
Therefore %
$$ \frac1v\le2\frac{m_k+n_k}{L_k}\toas{k\to\infty}2\den(x) ,$$
which proves the desired claim that $\hat{t}=\infty$ implies
$\den(x)\ge\frac1{2v}$.

Consider now the case of densities greater than $\frac1{2v}$. In
this case there might be two possibilities:

(a) All particles will eventually achieve the largest available
velocity $v$. Then the gaps will become not smaller than $v$ and
hence they cannot exceed $2v$ (by the assumption on the density
region). Obviously this situation may take place only if
$\den(x)\in[\frac1{2v},\frac1{v}]$ and it corresponds to the upper
branch of the Fundamental Diagram on Fig.~\ref{f:velocity-strong}.

(b) For any moment of time the are infinitely many particles
having gaps smaller than $v$ (and hence zero normalized local
velocities). Therefore at least for spatially periodic
configurations we can apply Lemma~\ref{l:max-gap} which guarantees
that only gaps smaller than $2v$ will survive with time. Thus to
study asymptotic properties it is enough to consider
configurations having only two types of gaps: smaller than $v$ and
between $v$ and $2v$.

Denote by $X(L,m,n)$ the subset of admissible configurations
$x\in X$ being spatially periodic with the spatial
period of length $L\in\IR_+$, which contains exactly $m\in\IZ_+$
particles with gaps belonging to the interval $[0,v)$ and
$n\in\IZ_+$ particles with gaps belonging to the interval $[v,2v)$.
Obviously $\den(x)=(m+n)/L$. The set $X(L,m,n)$ is invariant under
dynamics (each time when the size of a gap crosses the threshold $v$
one ``small'' gap becomes large and one ``large'' gap becomes
``small'') which immediately yields the exact value of the average
velocity $V(x)=\frac{nv}{m+n}$. On the other hand, by definition
$mv+n2v>L$ since the corresponding gaps fill in the segment of
length $l$ and lengths of both types of
gaps are smaller than $v$ and $2v$ respectively. Therefore %
$(\den(x)L+n)v>L$ and hence $n>L/v-\den(x)L$, which gives the
lower bound %
$$V(x) = \frac{nv}{m+n} = \frac{nv}{\den(x)L}
      > v~\frac{L/v-\den(x)L}{\den(x)L} = 1/\den(x) - v.$$ %
Observe, that choosing ``small'' and ``large'' gaps of length
$v-\ep$ and $2v-\ep$ for $0<\ep\ll1$ we see that the lower bound
can be ``almost'' achieved.

The upper bound of the average velocity in the hysteresis phase
(i.e. when $\frac1{2v}<\den(x)<\frac1{v}$) follows from the
existence of configurations with equal gaps of size larger than
$v$ for all densities from this segment. For the case
$\den(x)>1/v$ the upper bound is calculated using the opposite
length estimate $nv<L$. Then we get %
$$V(x) = \frac{nv}{m+n} = \frac{nv}{\den(x)L} < \frac{L}{\den(x)L}
       = 1/\den(x),$$ which agrees with the weak normalization case.

It remains to show that the region $H$ is filled in densely by the
pairs $(\den,V)$ corresponding to admissible configurations. To
this end one considers all possible choices of the integer
parameters $n,m$ and lengths of the corresponding gaps to get the
result. Indeed, $\forall \den\in(\frac1{2v},\frac1{v})$ there
exists an arbitrary large $L$ such that $\den L\in\IZ_+$. Choosing
now various available combinations of positive integers $m,n$ for
which $m+n=\den L$ we can approximate $V$ with the accuracy %
$$|V - \frac{nv}{m+n}|\le \frac{v}{\den L} \toas{L\to\infty}0.$$
\qed %

\n{\bf Remark}. By Theorem~\ref{t:fund-d-strong} for a.e. pair
$(V,\rho)\in$~H there exists an admissible configuration $x\in X$
such that $\den(x)=\rho$ and $V(x)=V$. On the other hand, it might
be possible that for some configurations having densities
belonging to the hysteresis region the average velocity is not
well defined and we claim only that all limit points of finite
time velocities belong to the vertical segment corresponding to
the given density.

\begin{corollary}\label{c:V-r-s} Let $x(r)\in X(r),~r>0$ and
$\den(x(r))$ be well defined and let ~$\forall i,t~v_i^t\equiv v$.
Then $V(x)=v$ if $\den(x(r))<\frac1{2v+2r}$ and otherwise
for a.e. point $(\den,V)$ in the curvilinear region %\newline%
$$H:=\{(\den,V): ~\max(\frac1{\den-2r}-v,0)\le V
     \le\min(\frac1{\den-2r},v)\}$$ %
there exists a configuration $x(r)\in X(r)$ with
$\den(x(r))=\den, V(x(r))=V$, i.e. the region $H$ corresponds to
the hysteresis.
\end{corollary}

\section{Local velocities of both signs}\label{s:vel-2signs} %
A close look to the previous analysis shows that we practically
did not use the property that all particles move in the same
direction, i.e. that $P(v_i^t\ge0)=1$. Now we explain the changes
necessary to study this more general case. Consider an infinite
configuration $x(r)\in X(r)$ and again interpret the values
$\{v_i^t\}_{i,t}$ (which now may have both positive and negative
signs, but still assuming that $|v_i^t|\le v$) as local velocities
for particles in the configuration $x^t(r)$.

The presence of particles moving in opposite directions leads to a
serious modification of the inequalities describing the violation
of the admissibility condition for the $i$-th local velocity.
Actually this is the main and the most serious change comparing to
the case of nonnegative velocities. Now we need to take into
account not only the position of the succeeding particle, but also
its velocity, as well as the corresponding quantities related to
the preceding particle. In this more general case the $i$-th local
velocity does not break the admissibility condition if and only if %
%
%\bea{\max(x^t_{i-1}(r), x^t_{i-1}(r)+v_{i-1}^t) + r %
%\a\le\min(x^t_i(r), x^t_i(r)+v_i^t) - r \\ %
%\a< \max(x^t_i(r), x^t_i(r)+v_i^t) + r %
%\le\min(x^t_{i+1}(r), x^t_{i+1}(r)+v_{i+1}^t) - r .}  %
%
\bea{\a\max(x^t_{i-1}(r), x^t_{i-1}(r)+v_{i-1}^t) + r  %
\le\min(x^t_i(r), x^t_i(r)+v_i^t) - r \\ %
\a\quad~< \max(x^t_i(r), x^t_i(r)+v_i^t) + r  %
\le\min(x^t_{i+1}(r), x^t_{i+1}(r)+v_{i+1}^t) - r .}  %
If for some $i\in\IZ$ and $j\in\{i-1,i+1\}$ the corresponding
inequality is not satisfied we say that there is a {\em conflict}
between the $i$-th particle and the $j$-th one and one needs to
resolve it. In terms of gaps $\G_i^t$ between particles the
inequalities above can be rewritten
as follows: %
\beq{e:adm-2signs}{\G_{j}^t\ge\max(v_j^t,~-v_{j+1}^t,~v_j^t-v_{j+1}^t),
                   ~~j\in\{i-1,i\} } %
Since the dynamics again will depend only on the sequence of gaps
$\{\G_i^t\}$ between particles, for each $r>0$ one can make the
invertible change of variables (\ref{e:r->r'}) (described in the
Introduction) to the case of `point' particles with $r=0$ which
we shall study further.

Exactly as in Section~\ref{s:intro} the {\em strong normalization}
means that we reject (nullify) all velocities leading to a
conflict, i.e %
$$ \cN_s(v_i^t,x^t):=\function{
      v_i^t &\mbox{if (\ref{e:adm-2signs}) holds } \\
      0 &\mbox{otherwise }.} $$
The situation with the weak normalization is more delicate. The
way how it was defined in Section~\ref{s:intro} can be
characterized as the only non-anticipating procedure allowing
conflicting particles to move simultaneously whenever possible.
Following this idea we say that a normalization is {\em weak} if
the positions of particles at the next time step
$x_i^{t+1}:=x_i^t+\cN_w(v_i^t,x^t)$ satisfy the conditions: %
\beq{e:weak-n}{ x_i^{t+1}\in\function{
      \{x_i^t+v_i^t\} &\mbox{if (\ref{e:adm-2signs}) holds } \\
      \{x_j^t, x_j^{t+1}\} &\mbox{if $\exists$ a conflict of the
                 particle $i$ with the particle $j=i\pm1$} .} }%
The 1st line describes the case when the admissibility condition
holds, while the 2nd line shows what happens if it breaks down.
Namely, if the $i$-th particle moves in the same direction as the
$j$-th one then (by the non-anticipation property) the former
assumes the previous position of the latter ($x_i^{t+1}=x_j^t$),
otherwise the positions of the conflicting particles at time $t+1$
coincide. The latter fact is the most important property here.

If directions of all instant local velocities coincide then
(\ref{e:weak-n}) defines the normalization uniquely. However if
their signs are different then (\ref{e:weak-n}) implies only that %
$$x_i^{t+1}=x_j^{t+1}\in[x_i^t,x_j^t]\cap[x_i^t+v_i^t,x_j^t+v_j^t].$$
Thus the set of weak normalizations is quite broad, for example it
includes a random normalization when two mutually conflicting
particles moving in opposite directions meet at a random point
belonging to the segments described above. One can give a
``natural'' specific construction of $\cNw$ normalizing local
velocities in such a way that positions of particles at the next
moment of time will be the same as if the particles would move
simultaneously at continuous time with the given local velocities
until the admissibility condition breaks
down: %
$$ \cN_{w,c}(v_i^t,x^t):=\function{
      v_i^t &\mbox{if (\ref{e:adm-2signs}) holds  }  \\
      -\G_{i-1}^t  &\mbox{if }
                   \G_{i-1}^t<-v_i^t,~ v_i^t<0,~v_{i-1}^t\le0 \\
      \G_{i}^t     &\mbox{if }
                   \G_{i}^t<v_{i}^t,~ v_i^t>0,~v_{i+1}^t\ge0 \\
      \frac{\G_{i-1}^t}{v_{i-1}^t-v_i^t}\times v_i^t  &\mbox{if }
                   \G_{i-1}^t<v_{i-1}^t-v_i^t,~ v_i^t<0,~v_{i-1}^t>0 \\
      \frac{\G_{i}^t}{v_{i}^t-v_{i+1}^t}\times v_i^t  &\mbox{if }
                   \G_{i}^t<v_{i}^t-v_{i+1}^t,~ v_i^t>0,~v_{i+1}^t<0 .}
$$

After this long discussion of the definition of the normalization
procedure it is surprising to find that all arguments used in the
analysis of the case of positive velocities remain valid with only
very slight changes.

\begin{lemma}\label{l:density-preservation+-}
The upper/lower densities $\den_\pm(x^t)$ are preserved under
dynamics.
\end{lemma}%
%  %
\proof One uses the same estimates as in the proof of
Lemma~\ref{l:density-preservation} except that now 2 particles may
simultaneously leave or enter a given spatial segment $I$ (instead
of 1). Thus the total change of the number of particles in
$I$ is less or equal to 2 and hence %
$$|\den(x^t,I) - \den(x^{t+1},I)|\cdot|I|\le2.$$ \qed%

\begin{lemma}\label{l:velocity-preservation+-}
Let $x\in X$ then $|V(x,j,t) - V(x,i,t)|\toas{t\to\infty}0$ a.s.
$\forall i,j\in\IZ$.
\end{lemma}%
\proof Again one follows the same argument as in the case of
nonnegative local velocities. The only difference is that in the
analysis of the connection between $\G_i^t$ and $\t\G_i^t$  now
one needs to consider new cases related to negative local
velocities.

Additionally here instead of the uniquely defined weak
normalization we need to consider an arbitrary one. If both
$v_i^t$ and $v_{i+1}^t$ are nonnegative we are in the situation
considered in Section~\ref{s:metric}. Therefore the cases (a) and
(b) hold automatically. Nevertheless we formulate all of them to
prove that $\t\G_i^t\ge\G_i^t~\forall t\in\IZ_0$: %
\begin{itemize}
\item[(a)] the condition (\ref{e:adm-2signs}) holds. Then
    obviously the argument used in Section~\ref{s:metric} woks.%
\item[(b)] $v_i^t>\G_i^t,~v_{i+1}^t\ge0$. Again one uses the same
    argument as in Section~\ref{s:metric}. %
\item[(c)] $v_i^t<-\G_{i-1}^t$. Then
    $\cNw(v_i^t,\t{x}^t)\le\cNw(v_i^t,x^t)\le0$ and
    $\cNw(v_{i+1}^t,\t{x}^t)\ge\cNw(v_{i+1}^t,x^t)$. Hence %
    \bea{\t\G_i^{t+1} \a= \t\G_i^{t} - \cNw(v_i^t,\t{x}^t)
                   + \cNw(v_{i+1}^t,\t{x}^t) \\ %
           \a\ge \G_i^{t} - \cNw(v_i^t,x^t) + \cNw(v_{i+1}^t,x^t) %
             =\G_i^{t+1}.} %
\item[(d)] $v_i^t\ge0,~v_{i+1}^t<0$ and $v_i^t-v_{i+1}^t>\G_i^t$.
  Then by definition $\t\G_i^{t+1}\ge0=\G_i^{t+1}$.
\end{itemize}

In the strong normalization setting one also considers the same
cases and proves by induction that $\t\G_i^{t+1}\ge\G_i^{t+1}-2v$
(instead of $\dots-v$ in the situation $v_i^t\ge0$). New cases are
the following %

\begin{itemize}
\item[(c')] $v_i^t<-\G_{i-1}^t$. Then
  $$\cNs(v_i^t,\t{x}^t)=v_i^t<-\G_{i-1}^t=\cNs(v_i^t,x^t)$$ and
  $$\cNs(v_{i+1}^t,\t{x}^t)-\cNs(v_{i+1}^t,x^t)\ge-2v$$ by the
  induction assumption. Hence %
  \bea{\t\G_i^{t+1} \a=\t\G_i^{t} - \cNs(v_i^t,\t{x}^t)
                      + \cNs(v_{i+1}^t,\t{x}^t) \\ %
      \a> \G_i^{t} - 2v  - \cNs(v_i^t,x^t) - \cNs(v_{i+1}^t,x^t) + 2v \\ %
      \a= \G_i^{t+1} .} %
\item[(d')] $v_i^t\ge0,~v_{i+1}^t<0$ and
  $\G_i^t<v_i^t-v_{i+1}^t\le\t\G_i^t$. Then %
  \bea{\t\G_i^{t+1}
  \a=\t\G_i^{t}+v_i^t-v_{i+1}^t>\t\G_i^{t}+\G_i^t\\
  \a\ge-2v+\G_i^t=-2v+\G_i^{t+1}.} %
\item[(d'')] $v_i^t\ge0,~v_{i+1}^t<0$ and
$\G_i^t<v_i^t-v_{i+1}^t>\t\G_i^t$. Then %
$$\t\G_i^{t+1}=\t\G_i^t\ge\G_i^t-2v=\G_i^{t+1}-2v.$$ %
\end{itemize}

Note that the difference $\t\G_i^{t+1}-\G_i^{t+1}=-2v$ may be
achieved only in the case (d').
The continuation of the proof is exactly the same as in
Section~\ref{s:metric}, except for the change of $2v$ to $4v$ in
the last inequality. \qed

Using these results and applying exactly the same arguments as in
the proof of Theorem~\ref{t:velocity-density} one gets the uniqueness
of the average velocity.

\begin{theorem}\label{t:velocity-density-2signs} %
In the weak normalization case the set of limit points as
$t\to\infty$ of the sequence $\{V(x,t)\}_{t\in\IZ_0}$ depends only
on the density $\den(x)$.
\end{theorem}%

\begin{theorem}\label{t:fund-d-weak-2signs} (Fundamental Diagram)
In the deterministic setting %
$V(x)=\lim\limits_{t\to\infty}\frac1t
         \sum\limits_{s=0}^{t-1}\min(1/\rho,v_0^s)$.
\end{theorem}
\proof Since at each moment of time $t\in\IZ_0$ the local
velocities of particles coincide, the condition (\ref{e:weak-n})
implies that
$$ x_i^{t+1}\in\{x_i^t+v_0^t,~x_{i\pm1}^t\} .$$
Thus the construction used in the proof of
Theorem~\ref{t:fund-d-weak} remains valid in this case as well.
\qed

\section{Generalizations and Discussion}\label{s:gen}

\subsection{Anticipating normalization} Throughout the paper we
consider only non-anticipating normalizations. In principle one
might try to consider an anticipating normalization allowing at
time $t$ the $i$-th particle to move up to the position of the
$(i+1)$-th particle $x_{i+1}^{t+1}$ at time $t+1$ rather than to
$x_{i+1}^t$. From the first sight this makes the normalization
scheme more flexible. Unfortunately the anticipating normalization
is not well posed since it turns out to be nonlocal. Namely a
single change in the sequence of local velocities (say of the
$i$-th one) may drastically alter the behavior of the system for
particles having indices arbitrary far from the changed one (i.e.
for $j\ll i$).

\subsection{One-sided particle densities}\label{s:one-sided} %
The density of a configuration in the way how it was defined in
Section~\ref{s:metric} depends sensitively on the statistics of
both left and right tails of the configuration. A close look shows
that in fact if all particles move in the same direction, say
right, one needs only the information about the corresponding
(right) tail, which allows to expand significantly the set of
configurations having densities and for which our results can be
applied.

For a configuration $x\in X$ by a {\em one-sided particle density}
we mean the limit %
\beq{e:one-side}
    {\hat\den(x):=\lim_{\ell\to\infty}\den(x,[0,\ell]).}%
The upper an lower one-sided densities correspond to the upper and
lower limits.

\begin{theorem}\label{t:one-side} Let $v_i^t\ge0~\forall i,t$.
Then all results of Lemma~\ref{l:den-r} and
Theorems~\ref{t:velocity-density}, \ref{t:fund-d-weak},
\ref{t:fund-d-strong} remain valid if one replaces the usual
particle density $\den$ to the one-sided density $\hat\den$.
\end{theorem}
\proof The key observation here is that the assumption
$v_i^t\ge0~\forall i,t$ implies that the movement of a given
particle in a configuration $x^t\in X$ depends only on particles
with larger indices. Therefore if one changes positions of all
particles with negative indices the particles with positive
indices will still have the same average velocity. On the other
hand, by Lemma~\ref{l:velocity-preservation} the average velocity
does not depend on the particle index. This allows to apply the
following trick.

For each configuration $x\in X$ of density $\den(x)$ we associate
a new configuration $\hat{x}\in X$ defined by the relation:
$$ \hat{x}_i:=\function{x_i^t &\mbox{if } i\ge0 \\
                        x_0+i/\den(x) &\mbox{otherwise }.} $$
Then obviously $\hat\den(x)=\den(\hat{x})=\den(x)$.

Therefore for all purposes related to the average velocities all
results valid for the configuration $\hat{x}$ remain valid for $x$
as well. \qed

Note however that this trick does not work for the case of local
velocities of both signs (considered in
Section~\ref{s:vel-2signs}), nor in the passive tracer analysis
(Section~\ref{s:tracer}). In both these situations statistics of
particles with negative indices cannot be neglected.

\subsection{Nagel-Schreckenberg traffic flow model} %
The celebrated Nagel-Schreckenberg traffic flow model introduced
in \cite{NS} for the lattice case is very similar to our case but
additionally to the lattice setting it uses a bit different
dynamics. In our terms this model differs from the main model
introduced in Section~\ref{s:intro} by that at each time step the
previous normalized local velocity of the $i$-th particle is
increasing by $0<a\le a_i^t$ until it reaches $v$. One can think
about $a_i^t$ as an acceleration under the action of a (random)
force (see e.g. \cite{Bl-hys}). Nevertheless the formalism
elaborated in the present paper allows to study the continuum
version of the Nagel-Schreckenberg model as well. In particular,
in the weak normalization case one applies basically the same
arguments as in Sections~\ref{s:metric},\ref{s:pairing} and
\ref{s:weak} since the distance between pair members cannot exceed
$C(v,a)\le v^2/a$. Note however that the average velocity should
be calculated in a more complicated way. Observe also that one can
consider random accelerations of both signs
$a_i^t\in(-\infty,-a]\cup[a,\infty)$ which makes the model more
applicable.

Mathematical formalism developed in the present paper can be
applied with minimal changes to a number of other traffic flow
models (discussed in detail, e.g. in a recent review \cite{MM})
allowing not only to study their continuum versions but also to
get rigorous results in the original lattice setting which are
absent at present.

\subsection{Passive tracer} \label{s:tracer}

Following the idea introduced in \cite{Bl-erg} we study the
dynamics of a passive tracer in the flow of particles imitating a
motion of a fast pedestrian in a slowly moving crowd of people.

Consider a pure deterministic setting ($v_i^t\equiv v$) with the
weak normalization and let $\map_v^tx$ describe the flow of
particles. The passive tracer occupies the position $y^t\in\IR$ at
time $t$ and moves all the time in the same direction. Before
carrying out the next time step of the model describing the flow
of particles, the tracer moves in its chosen direction to the
closest (in this direction) position of a particle of the
configuration $\map_v^tx$. After that the next iteration of the
flow occurs, the tracer moves to its new position, etc.

To be precise, let us fix a configuration $x\in X$ with
$\den(x)>0$ and introduce the maps
$\tau_{x}^{\pm}:\IR\to\IR$ defined as follows: %
$$\tau_{x}^{+}y := \min\{x_i: \; x_i>y\}, \quad %\CR
 \tau_{x}^{-}y := \max\{x_i: \; x_i<y\}.$$ Then the
simultaneous dynamics of the configuration of particles
(describing the flow) and the tracer is defined by the skew
product of two maps -- the map $\map_v$ and one of the maps
$\tau_{\cdot}^{\pm}$, i.e. %
$$(x,y) \to {\cal T}_{\pm}(x,y) := (\map_v x, \tau_{x}^{\pm}y),$$ %
acting on the extended phase space $X\times\IR$. The sign $+$ or
$-$ here corresponds to the motion along or against the flow. We
define the {\em average (in time) velocity} of the tracer
$$V_{{\rm tr}}(x,t):=(y^t-y^0)/t,$$ i.e. the total
distance covered by the tracer (which starts at position
$y^0\in\IR$) up to time $t\in\IZ_+$ with the positive sign if the
tracer moves forward, and the negative sign otherwise.

\begin{theorem} Let $v_i^t\equiv v~\forall i,t, ~~ \cN\equiv\cN_w,
~~ x\in X$ and let $x_{i+1}^0>x_{i}^0~~\forall i\in\IZ$.
If the tracer moves along the flow
(i.e. in the case ${\cal T}_+$) then
$$V_{{\rm tr}}(x,t)\toas{t\to\infty}V(x)=
  \function{v  &\mbox{if }~~ 0<\den(x)\le1/v\\
            1/\den(x) &\mbox{othewise } } .$$
If the tracer moves against the flow (case
${\cal T}_-$) then %
$V_{{\rm tr}}(x,t)\toas{t\to\infty}V(x)-1/\den(x)$.
\end{theorem}
\proof The assumption $x_{i+1}^0>x_{i}^0~~\forall i\in\IZ$ implies
that $x_{i+1}^t>x_{i}^t~~\forall i,t$ which allows to avoid a
pathology related to the presence of several particles at the same
position. In such a situation the tracer may ``jump'' through all
of them in one time step. This cannot happen if $r>0$ in
distinction to the case of point particles ($r=0$).

In the case of ${\cal T}_+$ the tracer will run down one of the
particles in the flow and will follow it, but cannot outstrip.
Thus $V_{{\rm tr}}(x,t)\toas{t\to\infty}V(x)$.

Consider now the case when the tracer moves backward with respect
to the flow, i.e. ${\cal T}_-$. Each time when the tracer
encounters a particle, on the next time step this particle moves
in the opposite direction and never will interfere with the
movement of the tracer. Thus during time $t>0$ the tracer meets
exactly $t$ particles which gives
$$(-V_{{\rm tr}}(x,t) + V(x,t))t\den(x) = t.$$
Therefore $$V_{{\rm tr}}(x,t)=-1/\den(x)+V(x,t).$$ %
\qed

Using similar arguments in the case of the strong normalization
one can show that $V_{{\rm tr}}(x,t)$ in the gaseous phase of the
particle flow has the same asymptotic as in the weak normalization
case. Since the flow in the fluid phase demonstrates hysteresis
the same phenomenon is unavoidable for the passive tracer as well.

\subsection{Multidimensional generalization}%
The constructions used in this paper are essentially
one-dimensional. Still at least some direct generalizations are
possible. Let $x_i^t\in\IR^d,~d\in\IZ_+$ and denote by $(x_i^t)_j$
the $j$-th coordinate of the $d$-dimensional vector $x_i^t$. We
say that a configuration $x^t(r)$ is {\em admissible} if %
\beq{e:admissible-mult}{\max_j((x_i^t(r))_j)+r\le
\min_j((x_{i+1}^t(r))_j)-r\qquad \forall i\in\IZ .}%

All results of Sections~~\ref{s:metric}, \ref{s:pairing},
\ref{s:vel-uniq}, and \ref{s:vel-2signs} hold in this setting.
Unfortunately the assumption~(\ref{e:admissible-mult}) implies
that a natural multidimensional generalization of the notion of
density of the configuration $x^t(r)$ turns out to be equal to
zero for any admissible configuration. However densities for
one-dimensional projections are well defined and for them the
Fundamental Diagram type results are readily available.

\subsection{Open problems and conjectures}

Our construction give a very precise information about the
asymptotic properties of DDS under consideration in the
deterministic setting. In the random setting we prove only the
uniqueness of the average velocity. From the results of
Section~\ref{s:metric} it follows that the mathematical
expectation of lower/upper average velocities are well defined but
we are not able to calculate them. On the other hand, we can
formulate a conjecture that the limits as time goes to infinity of
finite time average velocities are deterministic. In other words,
the Law of Large Numbers is valid for the sequence of finite time
average velocities.

An important question is whether the dynamical coupling of pairs
of processes with equal densities under the weak normalization is
nearly successful. Let $\cV$ be the common distribution of the
i.i.d. local velocities. As we know in the pure deterministic
setting when the distribution $\cV$ is concentrated at a single
point $\{v\}$ the dynamical coupling needs not to be successful.
Nevertheless we conjecture that for each nontrivial distribution
$\cV$ the nearly successful coupling takes place. Moreover, the
non-triviality of the distribution $\cV$ should lead to the
existence and uniqueness of the translationally invariant measure
of the Markov chain described by the DDS. Proofs of results of
this sort need the development of an additional probabilistic
apparatus and will be discussed elsewhere.

%\newpage%\small\footnotesize

\end{document}